\documentclass[review,12pt]{elsarticle}

\usepackage{amssymb}
\usepackage{mathtools} 
\usepackage[margin=2.54cm]{geometry} 
\usepackage{amsfonts} 
\usepackage[super]{nth} 
\usepackage{xcolor} 
\usepackage{booktabs} 
\usepackage{multirow} 
\usepackage{subcaption} 
\usepackage{tabularx} 
\usepackage{float} 
\usepackage{xr} 

\restylefloat{table} 
\restylefloat{figure} 

\journal{Energy}

\begin{document}

\setlength{\abovedisplayskip}{-12pt}
\setlength{\belowdisplayskip}{8pt}
\setlength{\abovedisplayshortskip}{-12pt}
\setlength{\belowdisplayshortskip}{8pt}

\begin{frontmatter}



\title{Optimal Control to Handle Variations in Biomass Feedstock Characteristics and Reactor In-Feed Rate}


\author[a]{Fikri Kucuksayacigil}
\author[a]{Mohammad Roni \fnref{*}}
\author[b]{Sandra D. Eksioglu}
\author[c]{Qiushi Chen}
\address[a]{Energy and Environment Science \& Technology Department, Idaho National Laboratory, Idaho Falls, ID}
\address[b]{Industrial Engineering Department, University of Arkansas, Fayetteville, AR}
\address[c]{Glenn Department of Civil Engineering, Clemson University, Clemson, SC}

\fntext[*]{Corresponding author: mohammad.roni@inl.gov}

\author{}


\address{}

\begin{abstract}

The variations in feedstock characteristics such as moisture and particle size distribution lead to an inconsistent flow of feedstock from the biomass pre-processing system to the reactor in-feed system. These inconsistencies result in low on-stream times at the reactor in-feed equipment. This research develops an optimal process control method for a biomass pre-processing system comprised of milling and densification operations to provide the consistent flow of feedstock to a reactor's throat. This method uses a mixed-integer optimization model to identify optimal bale sequencing, equipment in-feed rate, and buffer location and size in the biomass pre-processing system. This method, referred to as the hybrid process control (HPC), aims to maximize throughput over time. We compare HPC with a baseline feed forward process control. Our case study based on switchgrass finds that HPC reduces the variation of a reactor's feeding rate by 100\% without increasing the operating cost of the biomass pre-processing system for biomass with moisture ranging from 10 to 25\%. A biorefinery can adapt HPC to achieve its design capacity.

\end{abstract}




\end{frontmatter}


\section{Introduction} \label{sec:introduction}

Biomass characteristics vary by moisture, composition, particle size, shape, and density. These variabilities often cause poor flowability of materials at pre-processing systems that lead to clogging of equipment \cite{us2016biorefinery}. For example, variable moisture content causes variability in particle sizes. Whereas small particles may plug weep holes, larger particles may plug bins. Poor flowability often stems from arching, bridging, and rat-holing in the system \cite{us2016biorefinery}. Solid materials have unanticipated physical interactions with equipment, which often cause interruptions in the flow. These interruptions negatively affect down-steam processes.

The variability of biomass characteristics impacts biorefineries in different ways. Variations of particle size and distribution, moisture, ash, and heat contents negatively affect the integration of biomass feeding systems and conversion processes and result in low or unreliable on-stream time and long start-up times. For these reasons a number of biorefineries, which were planning to ramp-up production, had to delay these plans. For example, Bell \cite{bell2005challenges} and Merrow \cite{merrow2000developing} show that feeding systems reached only 77\% of their designed production amount at the end of the first year of operation (see \cite{bell2005challenges} for findings related to general particulate processes and \cite{merrow2000developing} for initial efforts to reveal and communicate drawbacks of particulate processes). Bioenergy experts report that capital expenditures and operating expenses have risen due to operational difficulties of handling biomass. These findings are proven empirically \cite{williams2016sources} and via engineering perspectives \cite{westover2018biomass}. The challenges outlined here emerge as significant barriers for biomass to become a sustainable source of energy. As a consequence, major commercial biorefineries (e.g., DuPont) have discontinued their production in the U.S. \cite{hirtzer2017}. Furthermore, investors are becoming reluctant to invest in the biofuel industry due to the high risks stemming from variability of biomass characteristics, market, and supply chain \cite{mamun2020supply}. Sustainable operations of a biorefinery under biomass variations remains a challenge \cite{hess2020advancements}.

A number of efforts have been made to mitigate the impacts of variations in biomass characteristics and to improve the design throughput of a biofuel conversion process. First, new designs of the equipment in the pre-possessing system have been proposed to improve flowability. There is rich literature examining physical interactions of biomass and equipment design to gain insight towards improving flowability. As an example, Dooley et al. \cite{dooley2013woody} proposed a rotary shear mechanism to chop biomass efficiently. Engineers alternatively prefer improving existing equipment, usually supported by modeling tools such as the discrete element method (DEM). For example, Ketterhagen et al. \cite{ketterhagen2009predicting} used DEM to identify the correct geometry of a hopper to ensure uninterrupted flow of particles. Crawford et al. \cite{crawford2016effects} investigated the effects of chemical additives, moisture content, particle interaction, and wall friction forces on the flowability of biomass materials. Wu et al. \cite{wu2011physical} measured physical characteristics of biomass and introduce designs of handling and storage equipment that align with these measurements.

Second, process control via advanced modeling and simulation tools has shown to reduce the negative impacts that variations of biomass characteristics have on achieving the design throughput of a conversion process. Process control shows a prominent potential to identify robust process designs that handle variations in feedstock characteristics and equipment down-time. There are abundant studies in the literature focusing on process control that model inflow and outflow for equipment in consideration. Marino et al. \cite{marino2017data, marino2018interpretable} developed a data-driven decision support systems, which inform operators regarding the performance of grinders given the moisture content of bales, in-feed rate, and screen size. Jaramillo and Sanchez \cite{jaramillo2018mass} proposed a mathematical formulation to model outflow from a pretreatment continuous tubular reactor. Numbi and Xia \cite{numbi2015systems} adopted a system optimization perspective to model optimal processing of high-pressure grinding rolls. Hartley et al. \cite{hartley2020effect} emphasized the effect of biomass physical characteristics on the performance of a pre-processing system by modeling biomass flow via a discrete-event simulation.

Third, taking preventative measures reduces equipment down-time due to variations of biomass characteristics. For instance, processing high-moisture bales at a lower rate minimizes the risk of clogging. We regard preventative measures as practical applications that improve the reliability of biomass pre-processing systems. One area in reliability is concerned with minimizing variability in production systems by adjusting control variables (e.g., feed at a lower rate). Strong mathematical assumptions, however, limit the applicability of the models in this area. Das et al. \cite{das2011volume} analyzed a production system where random failures of equipment occur. They provide a closed-form expression of the economic production lot-size. Li and Meerkov \cite{li2000production} suggested bounds on throughput variances so that a proper serial line of production can be designed. Assaf et al. \cite{assaf2014analytical} and Colledani et al. \cite{colledani2008analysis} analyzed the variability of output from unreliable equipment in single-stage or simplified multi-stage systems.

Last, but not least, changing the underlying structure of feeding systems by relocating equipment and storage areas has shown to improve the design throughput of a conversion process. For example, adding buffers to the system ensures a smooth flow. Indeed, improving production systems via investigating the effectiveness of buffer size and location is a rich area of research. A number of studies make probabilistic assumptions concerning processing time, biomass supply, etc. These models, however, are generally limited due to strong assumptions. Tempelmeier \cite{tempelmeier2003practical} investigated (continuous) flow of materials in production systems and discussed different algorithms to determine optimal buffer sizes and their locations. Tan and Gershwin \cite{tan2009analysis} analyzed a two-stage continuous-flow system by assuming that transitions among states are Markovian and one buffer with finite capacity is present between stages. Similarly, Turki et al.  \cite{turki2013perturbation} determined an optimal buffer size for a single-product machine with random failures. In Hosseini and Tan \cite{hosseini2017simulation}, each stage in a two-stage system has different states of processing with general distribution of transition times. The authors developed a mixed-integer linear programming model to determine optimal level of buffer by maximizing profit.

All efforts mentioned above have some limitations. Equipment design and process control are usually tied to a single equipment. In other words, these efforts tend to ignore the impacts of biomass flow to the performance of an integrated feeding and pre-processing system. Nevertheless, it is crucial that the design of a hopper or the control of a grinder are informed by the performance of predecessor and successor equipment. This is of significance to make system-wise optimal decisions. Moreover, studies concerning the reliability and buffer analysis are usually coupled with strong and unrealistic mathematical assumptions. Being tied to single equipment and strong mathematical assumptions makes the above-mentioned approaches inappropriate to model and analyze biomass feeding systems.

In our work, we mitigate the limitations of the above-mentioned approaches. This paper considers the whole biomass feeding system from beginning to end by tracking the flow dynamics and alterations in the physical characteristics of materials throughout the process. We develop and analyze process control strategies to handle the variation in feedstock characteristics and reactor in-feed rate, and do not restrict ourselves to strong modeling assumptions. As a result, our proposed model designs a reliable biomass feeding system for an integrated biorefinery.

We point out a few characteristics of our approach that distinguish it from other works in the literature. First, our approach considers sequencing of bales. Processing high-moisture bales often takes longer than usual due to the possibility of clogging in grinders. Such a practice becomes a bottleneck to reaching higher target rates. One way to overcome this challenge is to process low-moisture bales first to build-up inventory, and then, consume this inventory while processing high-moisture bales. Through this research, we analyze the impact that different sequences of bales have on achieving the maximum feeding rate, buffer size, and cost of operating the system. To the best of our knowledge, there is no study in the literature that have analyzed the effect of bale sequencing on the performance of feeding systems.

Second, our approach controls buffer size and location as the means for managing the variability of biomass characteristics and maintaining a continuous feeding of feedstock to a reactor at the maximum possible feeding rate. Finally, our work is concerned with the comparison of two different control approaches, the hybrid process control (HPC), which aims to maximize throughput over time, and a state-of-the-art baseline feed forward process control (BFFPC). Through this research, we provide cost/benefit analysis for {BFFPC} and {HPC} under different settings of a biomass feeding system. In summary, our study is uniquely positioned to devise a control approach for biomass feeding systems.

This research makes the following contributions to the literature. First, we propose {HPC} for biomass feeding systems that devises bale sequencing, as well as buffer location and size. We compare {HPC} with a state-of-the-art {BFFPC} for biomass pre-prossessing systems. Second, we exhibit how to formulate optimal control approaches using linear and mixed-integer linear programmings. These formulations are general to a large extent and can be applied to several biomass feeding systems. Third, we communicate a number of significant managerial insights. We convey that {HPC} is a potential solution to reduce variation of reactor feeding rate without increasing cost.

\section{Material and Methods} \label{sec:material_method}

\subsection{Feedstock Pre-processing Description}

\begin{figure}[H]
	\centering
	\begin{subfigure}[b]{\linewidth}
		\centering
		\includegraphics[width=\textwidth]{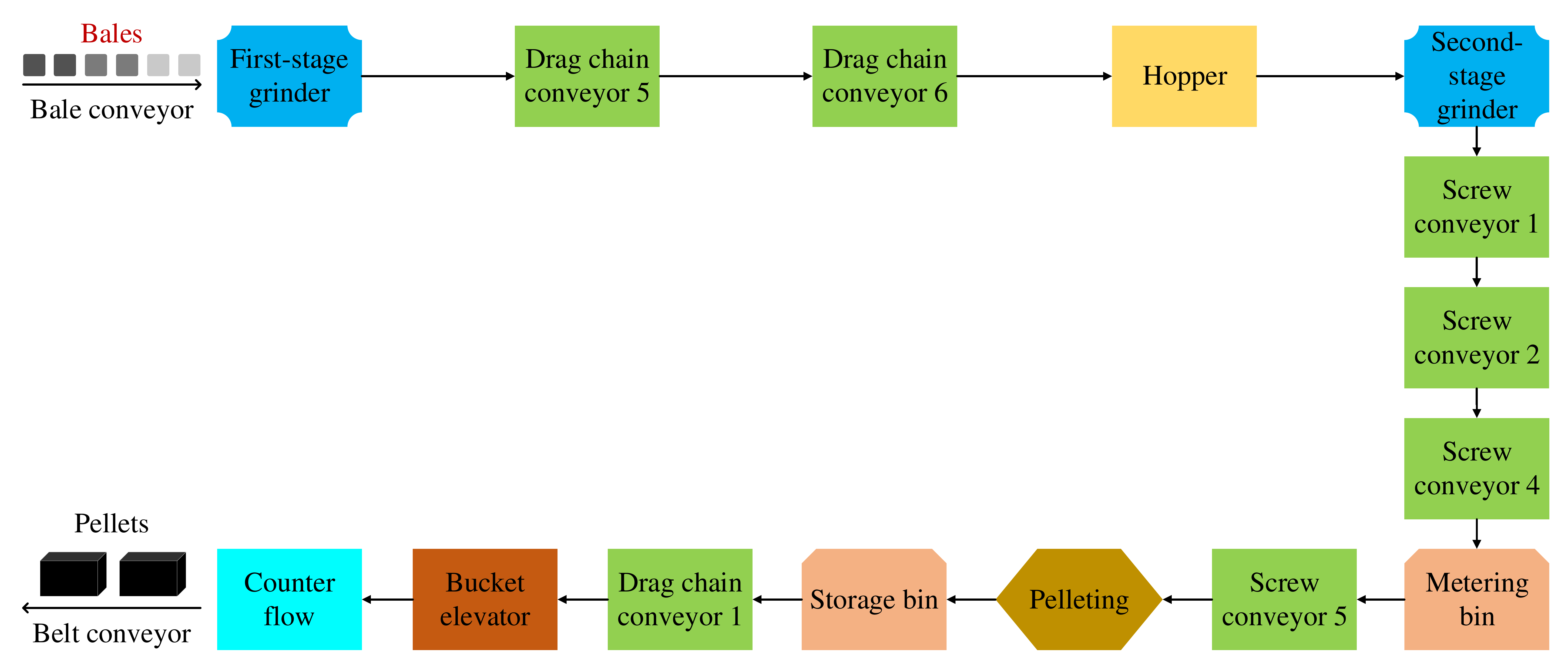}
		\caption{Without fractional milling}
		\label{fig:process_without}
		\vspace{2ex}
	\end{subfigure}
	\begin{subfigure}[b]{\linewidth}
		\centering
		\includegraphics[width=\linewidth]{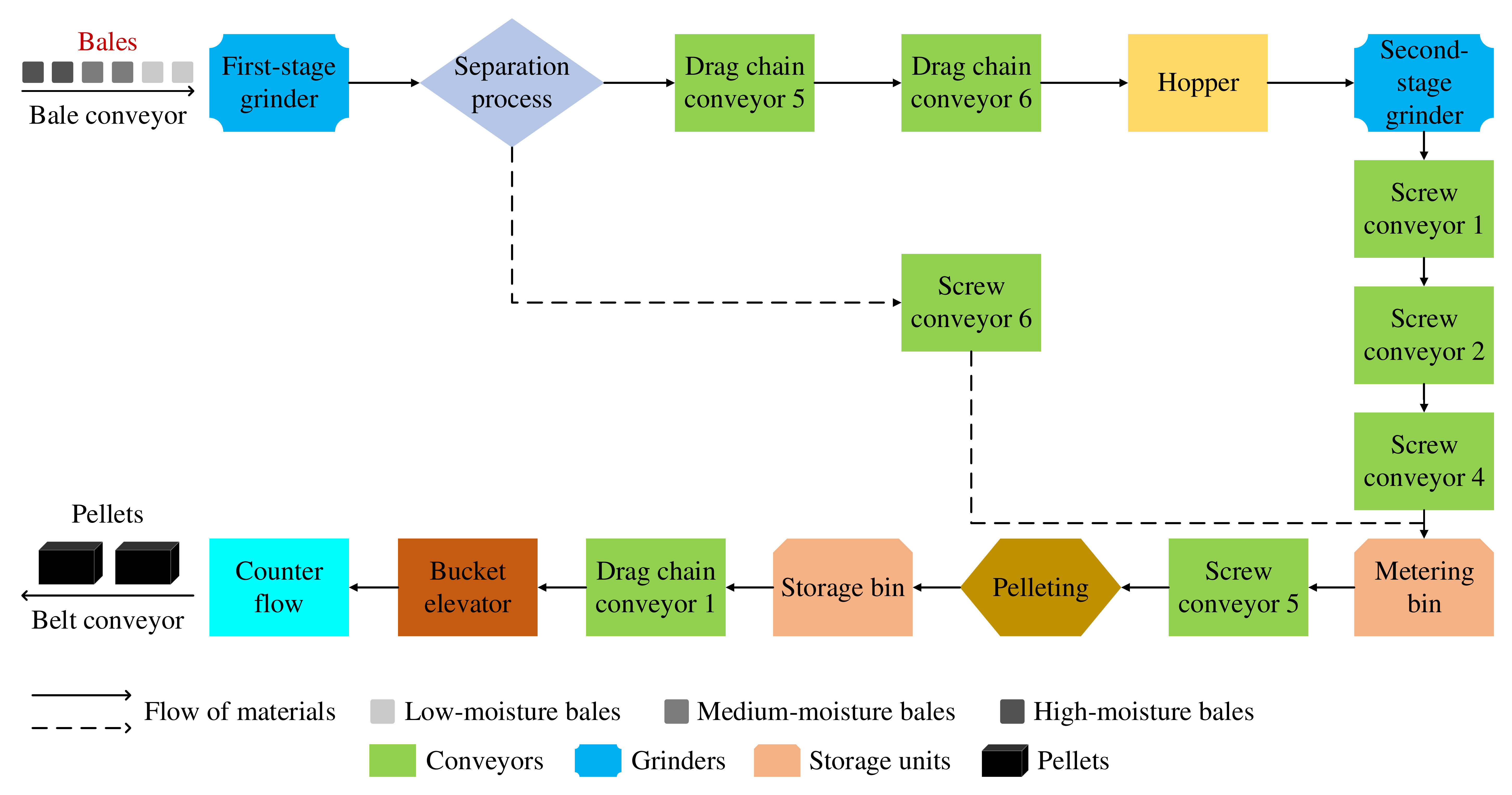}
		\caption{With fractional milling}
		\label{fig:process_with}
	\end{subfigure}
	\caption{A biomass pre-processing system}
	\label{fig:process}
\end{figure}

Feedstock pre-processing in this study includes two-stage size reduction and densification. The first-stage size reduction uses a grinder to break the bale into smaller particles. The second-stage size reduction further reduces the particle size via a hammer mill. Densification converts biomass to pellets. Figures \ref{fig:process_without} and \ref{fig:process_with} illustrate the system with equipment. We refer the interested readers to a study of biomass feedstock supply system design and analysis by \cite{under2014feedstock} for additional information about equipment and feedstock pre-processing systems. A supplemental file presents figures of several equipment and process demonstration unit (PDU). Table \ref{supp-tab:eqspecs} and \ref{supp-tab:costdata} list some specifications for equipment.

The feedstock pre-processing system that we utilize in this study has two different designs based on whether a separation unit is utilized between first-stage and second-stage grinding. Process without fractional milling does not have a separation process. After being processed in the first-stage grinder, materials go through the second-stage grinder for further size reduction. Such a practice produces fines that are prohibitive for the successful operation of the system. In process with fractional milling, following first-stage grinding, a separation unit separates small-sized particles and pushes them to a storage unit (we call the percentage of small-sized particles the ``bypass ratio"). Such a practice results in tighter particle size distribution, a minimum amount of fines, and reduced energy consumption \cite{kenney2013feedstock}.

\subsection{Description of Proposed Control Method in the Feedstock Pre-processing System}

Primary control variables in the feedstock pre-processing system are: (\textit{i}) screen sizes in grinders; (\textit{ii}) in-feed rate (amount of material processed per unit of time); (\textit{iii}) mill speed in grinders; (\textit{iv}) discharge rates from storage units; (\textit{v}) feeding rate of pellet mill; and (\textit{vi}) amount of moisture to add to pellet mix. Operators usually adjust these variables to ensure uniformity of particle sizes and increase pellet quality. Since screen size and moisture requirements can be determined in advance based on pellet particle size and quality requirements, we fix the value of the following control variables to predetermined values -- screen size associated with first- and second-stage size reduction and moisture.

The proposed control method (e.g., HPC) is an adaption of the feedback loop control and incorporates additional biomass process specific controls, examples of which are given in studies of refrigerators/cooling systems \cite{kyriakides2020dynamic} and demand-side management systems \cite{arnaudo2020techno}. Thus, HPC is a control approach whose primary goal is to meet a desire output response by adjusting bale sequencing decisions, biomass in-feed rate of different equipment, and buffer location and size. Consider the feedstock pre-processing design with fractional milling. Processing of bales with different moisture leads to additional operational constraints. For example, in-feed rate of second-stage grinding is constrained by moisture content of biomass since high moisture biomass is fed at a slow rate to prevent clogging in the system. Low moisture biomass is processed faster. This variation in moisture causes a fluctuation in the feedstock processing and reactor feeding rates.

The fluctuation of the reactor feeding rate can be handled in the proposed HPC method in the following way. First, the operator characterizes the bale based on its moisture range (e.g., high-, medium-, and low-moisture ranges). Second, the operator decides when to process low-, medium-, or high-moisture bales, as well as when to build-up/use the buffer to maintain a continuous feeding of the reactor. By simultaneous optimization of bale sequencing, in-feed rate, and buffer location and size, the HPC determines an optimal trajectory that maintains the continuous feeding of a reactor. If the current process buffer size (e.g., metering bin) is not sufficient, then optimization also decides by how much to increase the buffer size (see Figure \ref{supp-fig:hpc_control}).

The proposed HPC method is compared with the state-of-the-art {BFFPC} to demonstrate the differences and benefits of the proposed method. BFFPC faces the same operational constraints caused by moisture. This BFFPC method maximizes the total throughput over time but does not necessarily maintain a continuous feeding of the reactor. The BFFPC operates as follows. First, the operator characterizes bales based on moisture range, but feeds bales in an unsystematic way. That means, bales are not arranged in a sequence according to their moisture contents. We will henceforth use the term ``random {pattern}" to refer to this arrangement. Second, in contrast to the HPC method, BFFPC finds an optimal trajectory that maximizes throughput by optimizing in-feed rate, but it does not control buffer location and size (see Figure \ref{supp-fig:bffpc_control}).

We summarize the commonalities and differences of BFFPC and HPC in Table \ref{supp-tab:control_summary}.

\subsection{Optimization Model Development} \label{sec:optimization}

The optimization model identifies the minimum time to process a {fixed} number of bales using HPC and BFFPC under various operational constraints. The optimal trajectory to process a fixed number of bales is found via an {algorithm} as shown in Figure \ref{supp-fig:flowchart}. The goal is to shorten the modeling horizon iteratively until no better solution is found. In each iteration of this {algorithm}, an optimization model is solved. The implementation of this {algorithm} is described in supplementary material.

\subsubsection{MILP model for HPC} \label{sec:hpcformulation}

We use the following notation to model HPC. Throughout the text, we use dry Mega gram (dry Mg) as mass unit.

\underline{\textit{Sets:}}
\vspace{-0.5em}
\begin{itemize}
	\itemsep-0.5em
	\item $\mathcal{I}$: Equipment used in a biomass feeding system
	\item $\mathcal{I_{G}}$: Grinders
	\item $\mathcal{I_{M}}$: Storage units
	\item $\mathcal{I}_{i}$: Predecessor of equipment $i$
	\item $\mathcal{T}$: Periods in modeling horizon
	\item $\mathcal{L}$: Biomass moisture levels considered in model
	\item $\mathcal{C}$: Potential expansion options for storage units, $\mathcal{C} \coloneqq \{0, 0.1, 0.2, \ldots, 1\}$ where 0.1 means capacity is expanded by 10\%.
\end{itemize}

\underline{\textit{Indices:}}
\vspace{-0.5em}
\begin{itemize}
	\itemsep-0.5em
	\item $i \in \mathcal{I}$: Equipment locations in the system
	\item $t \in \mathcal{T}$: Periods
	\item $s \in \mathcal{L}$: Moisture level
	\item $k \in \mathcal{C}$: Potential expansion option.
\end{itemize}

\underline{\textit{Parameters:}}
\vspace{-0.5em}
\begin{itemize}
	\itemsep-0.5em
	\item $\alpha$: {Penalty} for adjustment of reactor feeding (\$/dry Mg/unit of time) 
	\item $\eta$: Unit sale price of pellets (\$/dry Mg)
	\item $c_{is}$: Unit cost of equipment $i$ based on moisture level $s$ (\$/unit of time)
	\item $c_{isk}$: Cost of operating equipment $i \in \mathcal{I_{M}}$ with processing capacity $k$ for biomass with moisture level $s$  (\$)
	\item $w$, $h$: Width and height of a bale (meters)
	\item $j_{st}$: An indicator parameter which takes the value 1 if biomass with moisture level $s$ is processed in period $t$, and takes the value 0 otherwise
	\item $d_{s}$: Density of biomass bale with moisture level of $s$ (dry Mg/cubic meter)
	\item $\bar{d}_{is}$: Average density of biomass with moisture level $s$ in equipment $i \in \mathcal{I_{M}}$ (dry Mg/cubic meter)
	\item $\overline{x}_{is}$: Capacity of equipment $i$ when processing materials with moisture level of $s$ (dry Mg/unit of time)
	\item $n_{s}$: The number of bales with moisture level $s$ 
	\item $q_{s}$: The mass of a bale with moisture level $s$ (dry Mg)
	\item $\mu_{i}$: Dry matter loss in equipment $i \in \mathcal{I_{G}}$ (\%)
	\item $\theta_{s}$: Bypass ratio when processing biomass with moisture level of $s$ (\%)
	\item $\overline{m}_{i}$: Mass capacity of equipment $i \in \mathcal{I_{M}}$ (dry Mg)
	\item $\overline{v}_{i}$: Volume capacity of equipment $i \in \mathcal{I_{M}}$ (cubic meter)
	\item $\gamma_{t} \triangleq w h \sum_{s \in \mathcal{L}} d_{s} j_{st}$: Amount of biomass per meter of bale conveyor at period $t$ (dry Mg/meter).
\end{itemize}

\underline{\textit{Decision variables:}}

\emph{Binary variables}
\vspace{-0.5em}
\begin{itemize}
	\itemsep-0.5em
	\item $J_{k}$: takes the value 1 if the capacity of an equipment is expanded by $k$\%, takes the value 0 otherwise.
\end{itemize}

\vspace{-0.5em}
\emph{Continuous variables}
\vspace{-0.5em}
\begin{itemize}
	\itemsep-0.5em
	\item $X_{t}$: In-feed rate of the system at period $t$ (dry Mg/unit of time)
	\item $X_{it}$: Outflow from equipment $i$ at period $t$ (dry Mg/unit of time)
	\item $\mathbb{X}_{ist}$: Outflow of material with moisture level $s$ from equipment $i \in \mathcal{I_{M}}$ at period $t$ (dry Mg/unit of time)
	\item $\mathbb{M}_{ist}$: Inventory level of biomass with moisture level $s$ in equipment $i \in \mathcal{I_{M}}$ at period $t$ (dry Mg)
	\item $Y_{t}$: {Conveyor speed} at period $t$ (meters/unit of time)
	\item $\beta_{t}^{+}$, $\beta_{t}^{-}$: Decision variables to control reactor feeding rate deviations (dry Mg/unit of time).
\end{itemize}

The optimization model we make reference to in Figure \ref{supp-fig:flowchart} is defined by Equations \eqref{eq:hybridobjective}-\eqref{eq:betaconst} and corresponds to HPC. The objective function \eqref{eq:hybridobjective} minimizes the variations of reactor's feeding, maximizes the throughput of the system, and determines optimal location and size of the buffer. The first term in the objective maximizes the throughput of the system by maximizing $\eta X_{|\mathcal{I}|t}$, the hypothetical revenue gained by producing/selling one unit of pellet at period $t$. The second term in the objective minimizes the volatility of the feeding rate of the reactor to ensure a continuous and smooth feeding of the reactor, a key success indicator for a biomass feeding system. The last term in the objective minimizes the cost of expanding the capacity of storage units.

\begin{align}
\max \sum_{t \in \mathcal{T}} \left[ \eta X_{|\mathcal{I}|t} - \alpha \left( \beta_{t}^{+} + \beta_{t}^{-} \right) \right] - \sum_{i \in \mathcal{I_{M}}} \sum_{k \in \mathcal{C}} \sum_{s \in \mathcal{L}} \sum_{t \in \mathcal{T}} c_{isk} j_{st} J_{k} \label{eq:hybridobjective}
\end{align}

We use Equations \eqref{eq:hybridobjective}-\eqref{eq:betaconst} to obtain optimal trajectories of system in-feed rate, reactor feeding rate, and inventory in storage units. Operating cost of equipment however has unit of \$/hour. Thus, after obtaining the optimal trajectories, we calculate the corresponding cost of operating the system via 

\begin{align}
\sum_{i \in \left( \mathcal{I} - \mathcal{I_{M}} \right)} \sum_{s \in \mathcal{L}} \sum_{t \in \mathcal{T}} c_{is} + \sum_{i \in \mathcal{I_{M}}} \sum_{k \in \mathcal{C}} \sum_{s \in \mathcal{L}} \sum_{t \in \mathcal{T}} c_{isk} j_{st} J_{k} \label{eq:costcalculation}
\end{align}

\noindent given the optimal value of $J_{k}$. Therefore, the optimal cost we present in our numerical analyses considers only system operating costs.

\textbf{Capacity constraints}: Each equipment has processing capacities, which we model via

\begin{align}
X_{it} & \leq \sum_{s \in \mathcal{L}} j_{st} \overline{x}_{is} & \forall i \in \mathcal{I}, \forall t \in \mathcal{T} \label{eq:outflowmax} \\
\sum_{s \in \mathcal{L}} \mathbb{M}_{ist} & \leq \overline{m}_{i} \left(1 + \sum_{k \in \mathcal{C}} kJ_{k} \right) & \forall i \in \mathcal{I_{M}}, \forall t \in \mathcal{T} \label{eq:masscapacityexpanded} \\
\sum_{s \in \mathcal{L}} \frac{\mathbb{M}_{ist}}{\bar{d}_{is}} &  \leq \overline{v}_{i} \left(1 + \sum_{k \in \mathcal{C}} kJ_{k} \right) & \forall i \in \mathcal{I_{M}}, \forall t \in \mathcal{T} \label{eq:volumecapacityexpanded} \\
1 & \geq \sum_{k \in \mathcal{C}} J_{k}. \label{eq:sumtoone}
\end{align}

Constraint \eqref{eq:outflowmax} indicate that the outflow of biomass from each equipment should not exceed its pre-determined processing capacity, which is driven by the moisture level of the biomass processed at each period. Constraint \eqref{eq:masscapacityexpanded} limit the mass of biomass inventoried and constraint \eqref{eq:volumecapacityexpanded} limit the volume of biomass inventoried in a storage unit. This consideration is necessary because for loose biomass (low density), the volume capacity is restrictive, whereas mass capacity turns out to be restrictive in the case when biomass is highly dense. These constraints also consider the impact that extensions of the storage capacity have on the mass and volume of the inventory.  We calculate the average density of biomass $\bar{d}_{is}$ as the weighted average of biomass densities conveyed via the two conveyors connected to a metering bin. The corresponding weights are $\theta_{s}$ and $1 - \theta_{s}$. Constraint \eqref{eq:sumtoone} ensures that only one capacity expansion is chosen. If we desire to run HPC without the capacity expansion option, we simply adjust constraint \eqref{eq:sumtoone} by replacing its left-hand side with 0.

\textbf{Operational constraints:} The system under consideration has operational capacities related to system in-feed. The in-feed rate of the system refers to the amount of biomass fed to the system in a period. Constraint \eqref{eq:infeedbeltspeed} calculates the in-feed rate based on the physical dimensions of a bale, its density, and {conveyor speed}. Constraint \eqref{eq:allprocessed} ensures that every bale is processed. $\beta_{t}^{+}$ and $\beta_{t}^{-}$ calculate the change of reactor feeding rate from one period to the next. If $X_{|\mathcal{I}|t} - X_{|\mathcal{I}|t-1} > 0$, $\beta_{t}^{+} > 0$ and $\beta_{t}^{-} = 0$, otherwise, if $X_{|\mathcal{I}|t} - X_{|\mathcal{I}|t-1} < 0$, $\beta_{t}^{+} = 0$ and $\beta_{t}^{-} > 0$. In the case when the feeding of reactor does not change in successive time periods, $\beta_{t}^{+} = \beta_{t}^{-} = 0$.

\begin{align}
X_{t} & = \gamma_{t} Y_{t} & \forall t \in \mathcal{T} \label{eq:infeedbeltspeed} \\
\sum_{t \in \mathcal{T}} j_{st} X_{t} & = q_{s} n_{s} & \forall s \in \mathcal{L} \label{eq:allprocessed} \\
\beta_{t}^{+} - \beta_{t}^{-} & = X_{|\mathcal{I}|t} - X_{|\mathcal{I}|t-1} & \forall t \in \mathcal{T} \label{eq:betacalc}
\end{align}

\textbf{Flow balance constraints:} Feeding system is essentially a flow system that needs to satisfy flow balance constraints.  Constraint \eqref{eq:flowgrinders} computes the outflow from the grinders. The outflow is reduced by the dry matter losses in the grinders. Constraint \eqref{eq:flowdcc5} shows that the amount of biomass that does not bypass will continue flowing on drag chain conveyor 5 (d.c.5). The biomass that bypasses the separation unit flows on screw conveyor 6 (s.c.6) that feeds the metering bin, as indicated by constraint \eqref{eq:flowsc6}. Constraint \eqref{eq:flowdrags} represents that flow balance (flow in equals flow out) for other equipment.

\begin{align}
X_{it} & = \sum_{p \in \mathcal{I}_{i}} (1 - \mu_{i}) X_{pt} & \forall i \in \mathcal{I_{G}}, \forall t \in \mathcal{T} \label{eq:flowgrinders} \\
X_{it} & = \sum_{p \in \mathcal{I}_{i}} \left( 1 - \sum_{s \in \mathcal{L}} j_{st} \theta_{s} \right) X_{pt} & i = \text{d.c.5}, \forall t \in \mathcal{T} \label{eq:flowdcc5} \\
X_{it} & = \sum_{p \in \mathcal{I}_{i}} \left( \sum_{s \in \mathcal{L}} j_{st} \theta_{s} \right) X_{pt} & i = \text{s.c.6}, \forall t \in \mathcal{T} \label{eq:flowsc6} \\
X_{it} & = \sum_{p \in \mathcal{I}_{i}} X_{pt} & \forall i \in \mathcal{I} - \{ \text{d.c.5, s.c.6}, \mathcal{I_{M}}, \mathcal{I_{G}} \}, \forall t \in \mathcal{T} \label{eq:flowdrags}
\end{align}

\textbf{Inventory balance constraints:} Constraint \eqref{eq:invall} is the mass balance constraint for storage units. Constraint \eqref{eq:flowmetering} finds the total outflow from each storage unit at period $t$. We assume that the system does not carry initial inventory in storage units, presented by constraint \eqref{eq:initialzero}.

\begin{align}
\mathbb{M}_{ist} & = \mathbb{M}_{ist-1} + \sum_{p \in \mathcal{I}_{i}} j_{st}X_{pt} - \mathbb{X}_{ist} & \forall i \in \mathcal{I_{M}}, \forall s \in \mathcal{L}, \forall t \in \mathcal{T} \label{eq:invall} \\
X_{it} & = \sum_{s \in \mathcal{L}} \mathbb{X}_{ist} & \forall i \in \mathcal{I_{M}}, \forall t \in \mathcal{T} \label{eq:flowmetering} \\
\mathbb{M}_{is0} & = 0 & \forall i \in \mathcal{I_{M}}, \forall s \in \mathcal{L} \label{eq:initialzero}
\end{align}

\textbf{Non-negativity constraints}: The following are the non-negativity constraints.

\begin{align}
0 & \leq X_{t} & \forall t \in \mathcal{T} \label{eq:xtnonnegative} \\
0 & \leq X_{it} & \forall i \in \mathcal{I}, \forall t \in \mathcal{T} \\
0 & \leq Y_{t} & \forall t \in \mathcal{T} \\
0 & \leq \mathbb{X}_{ist}, \mathbb{M}_{ist} & \forall i \in \mathcal{I_{M}}, \forall s \in \mathcal{L}, \forall t \in \mathcal{T} \label{eq:xmnonnegative} \\
0 & \leq \beta_{t}^{+}, \beta_{t}^{-} & \forall t \in \mathcal{T} \label{eq:betaconst}
\end{align}

We refer to the MILP formulation above (\eqref{eq:hybridobjective}, \eqref{eq:outflowmax} to \eqref{eq:betaconst}) as \textbf{Problem H} for HPC.

Notice that this version of HPC includes fractional milling. If we desire to use HPC without fractional milling, we make the following adjustments in the above model. We set $\theta_{s} = 0, \forall s \in \mathcal{L}$, which impacts constraints \eqref{eq:flowdcc5} and \eqref{eq:flowsc6}. Fixing the value of $\theta_{s} = 0$ also impacts constraint \eqref{eq:volumecapacityexpanded} since $\theta_{s}$ and $ 1 - \theta_{s}$ are used as weights in calculations of $\bar{d}_{is}$. Since screw conveyor 6 does not exist in this problem setting, we deduct the cost of screw conveyor 6 from the cost expression \eqref{eq:costcalculation}.

\subsubsection{LP model for BFFPC} \label{sec:bffpcformulation}

We develop an LP formulation of BFFPC by updating the model introduced above for HPC. In this section, we only list the modifications we make.

\textbf{Objective function}: The objective function of BFFPC focuses on maximizing the amount of biomass processed within the time horizon.

\begin{align}
\max \sum_{t \in \mathcal{T}} \left[ \eta X_{|\mathcal{I}|t} \right] \label{eq:baseobjective}
\end{align}

We compute the cost of operating the system via

\begin{align}
\sum_{i \in \mathcal{I}} \sum_{s \in \mathcal{L}} \sum_{t \in \mathcal{T}} c_{is} \label{eq:costofbffpc}
\end{align}

\textbf{Capacity constraints of storage units}: Since BFFPC does not consider expansion of  storage units capacities, we replace constraints \eqref{eq:masscapacityexpanded} and \eqref{eq:volumecapacityexpanded} of HPC with constraints \eqref{eq:masscapacity} and \eqref{eq:volume_metering}, as shown below.

\begin{align}
\sum_{s \in \mathcal{L}} \mathbb{M}_{ist} & \leq \overline{m}_{i} & \forall i \in \mathcal{I_{M}}, \forall t \in \mathcal{T} \label{eq:masscapacity} \\
\sum_{s \in \mathcal{L}} \frac{\mathbb{M}_{ist}}{\bar{d}_{is}} & \leq \overline{v}_{i} & \forall i \in \mathcal{I_{M}}, \forall t \in \mathcal{T} \label{eq:volume_metering}
\end{align}

Notice that BFFPC does not use variables $J_{k}$; thus, constraint \eqref{eq:sumtoone} is not considered. Since BFFPC does not aim to control the variations of a reactor's feeding rate, variables $\beta_{t}^{+}$, $\beta_{t}^{-}$ and constraints \eqref{eq:betacalc} and \eqref{eq:betaconst} are not considered as well.

The remaining of the constraints are the same in both HPC and BFFPC model formulations. We refer to the LP formulation  presented above (e.g., \eqref{eq:baseobjective}, \eqref{eq:outflowmax}, \eqref{eq:infeedbeltspeed}, \eqref{eq:allprocessed}, \eqref{eq:flowgrinders}-\eqref{eq:xmnonnegative}, \eqref{eq:masscapacity}, and \eqref{eq:volume_metering}) as  \textbf{Problem B} for BFFPC.

\subsection{Solution Approach}

We implement the LP and MILP models in Julia \cite{bezanson2017julia}, utilizing JuMP modeling language \cite{DunningHuchetteLubin2017}, and use the Gurobi solver \cite{gurobi2020}. We execute the models on a computer that has 32 GB RAM and an Intel(R) Xeon(R) CPU E5-1620 v2 CPU running at 3.7 GHz. We use a 1-minute granularity in the optimization. We set a 30-minute time limit and 1\% optimality gap for Gurobi. The minimum, average, and maximum CPU time to solve a single LP (\textbf{Problem B}) are 0.282, 0.405, and 0.642 seconds. The corresponding minimum, average, and maximum CPU time to solve the MILP \textbf{Problem H} are 0.318, 0.721, and 1.774 seconds. The running time of the algorithm presented in Figure \ref{supp-fig:flowchart} varies from 18.87 seconds to 66.59 seconds. This variation results from different patterns of feeding bales to the system.

\section{Case Study} \label{sec:casestudy}

\subsection{Technical Process Data}

We utilize technical process data from the PDU at Idaho National Laboratory (INL). Table \ref{tab:inputdata} summarizes the data collected by analyzing historical grinding and densification experiments of switchgrass. The data is grouped based on the moisture level of biomass and based on processes. Notice that grinders consume more energy when biomass has a high moisture level \cite{mani2004grinding}. Due to grinding, biomass loses moisture. Wet material loss is the largest for biomass with high moisture levels \cite{tumuluru2019biomass}. Since particle sizes alter during grinding, material density changes. Bulk density of biomass diminishes during first-stage grinding as biomass in bale form is tapped. In contrast, density increases during the second-stage grinding as the particles get smaller. As noted before, wet biomass is processed slow to prevent clogging. Hence, the maximum in-feed rate of grinders decreases for high-moisture materials, which negatively impacts the processing capacity of the system. There are several options to select screen size combinations for coupled grinders. Experiments at PDU show that it is the best practice to use screens with a diameter of 76.2-mm and 11.11-mm in first-stage and second-stage grinders, respectively \cite{tumuluru2017}. Utilizing these screens minimizes the amount of fine particles.

\begin{table}[H]
	\footnotesize
	\centering
	\begin{tabular}{lllccc}
		\toprule
		& & & \multicolumn{3}{c}{\textbf{Moisture Levels}} \\
		& & & \textbf{High} & \textbf{Medium} & \textbf{Low} \\
		\cline{4-6}
		\multirow{2}{*}{\textbf{Bale}} & \multicolumn{2}{l}{Moisture (\%)} & 25 & 17.50 & 10 \\
		& \multicolumn{2}{l}{Dry bulk density (dry Mg/cubic meter)} & 0.144 & 0.144 & 0.144 \\
		\midrule
		\multirow[b]{4}{*}{\textbf{First-stage}} & \multirow{3}{*}{Operating conditions} & Screen size (mm) & 76.2 & 76.2 & 76.2 \\
		& & Dry bulk density (dry Mg/cubic meter) & 0.144 & 0.144 & 0.144 \\
		& & Moisture (\%) & 25 & 17.5 & 10 \\
		\cline{2-6}
		& \multirow{5}{*}{Process performance} & Energy consumption (kWh/dry Mg) & 25.68 & 21.95 & 7.10 \\
		\multirow[t]{4}{*}{\textbf{Grinder}} & & Moisture loss (\%) & 4.77 & 3.00 & 0.5 \\
		& & Dry matter loss (\%) & 1.50 & 1.50 & 1.50 \\
		& & Bulk density change (dry Mg/cubic meter) & -0.09 & -0.10 & -0.11 \\
		& & Maximum in-feed rate (dry Mg/hour) & 2.20 & 4.53 & 5.23 \\
		\midrule
		\multirow[b]{3}{*}{\textbf{Separation}} & \multirow{3}{*}{Operating conditions} & Screen size (mm) & 6.35 & 6.35 & 6.35 \\
		& & Dry bulk density (dry Mg/cubic meter) & 0.053 & 0.041 & 0.039 \\
		& & Moisture (\%) & 20.23 & 14.50 & 9.50 \\
		\cline{2-6}
		\multirow[t]{2}{*}{\textbf{Unit}} & \multirow{2}{*}{Process performance} & Moisture loss (\%) & 0.71 & 0.50 & 0.00 \\
		& & Bypass ratio (\%) & 40.48 & 44.98 & 49.98 \\
		\midrule
		\multirow[b]{4}{*}{\textbf{Second-stage}} & \multirow{3}{*}{Operating conditions} & Screen size (mm) & 11.11 & 11.11 & 11.11 \\
		& & Dry bulk density (dry Mg/cubic meter) & 0.053 & 0.041 & 0.039 \\
		& & Moisture (\%) & 19.52 & 13.99 & 9.50 \\
		\cline{2-6}
		& \multirow{5}{*}{Process performance} & Energy consumption (kWh/dry Mg) & 49.71 & 16.58 & 14.67 \\
		\multirow[t]{4}{*}{\textbf{Grinder}} & & Moisture loss (\%) & 4.00 & 3.00 & 0.7 \\
		& & Dry matter loss (\%) & 0.50 & 0.50 & 0.50 \\
		& & Bulk density change (dry Mg/cubic meter) & 0.066 & 0.082 & 0.090 \\
		& & Maximum in-feed rate (dry Mg/hour) & 1.59 & 2.80 & 5.23 \\
		\midrule
		\multirow[b]{3}{*}{\textbf{Pelleting}} & \multirow{2}{*}{Operating conditions} & Dry bulk density (dry Mg/cubic meter) & 0.119 & 0.123 & 0.129 \\
		& & Moisture (\%) & 15.52 & 10.99 & 8.80 \\
		\cline{2-6}
		& \multirow{4}{*}{Process performance} & Energy consumption (kWh/dry Mg) & 90.39 & 60.63 & 55.12 \\
		\multirow[t]{3}{*}{\textbf{Equipment}} & & Moisture loss (\%) & 3.92 & 1.50 & 0.00 \\
		& & Bulk density change (dry Mg/cubic meter) & 0.547 & 0.542 & 0.537 \\
		& & Maximum in-feed rate (dry Mg/hour) & 3.33 & 3.81 & 4.76 \\
		\bottomrule
	\end{tabular}
	\caption{Technical process-related data used as input for BFFPC and HPC}
	\label{tab:inputdata}
\end{table}

Table \ref{tab:inputdata} also presents the technical process data for the separation unit. The bypass ratio is smallest in the case when biomass is highly moist. When materials are wet, they tend to be cohesive, cling to one another, and end up creating aggregates. This reduces the amount of materials bypassing the separation unit. Materials lose moisture at this stage as well. Wet materials lose the largest amounts of moisture during separation. Table \ref{tab:inputdata} also lists the technical process data for densification equipment (via pelleting mechanism). Pelleting equipment consumes the highest amount of energy in a biomass feeding system \cite{yancey2013drying}. When material is wet, energy consumption rises. Moisture loss occurs due to heating during densification. Wet materials lose the highest amounts of moisture during pelletization. Operators feed pelleting equipment at a lower rate when material is highly moist to ensure the production of high-quality pellets.

All equipment have a processing capacity of 13.61 dry Mg/hour with the exception of the first-stage grinder, second-stage grinder, and pelleting. Metering and storage bins have a storage capacity of 4.54 dry Mg. The volume capacity of these storage units is 24.07 cubic meter. Width, height, and length of a bale are 1.22, 0.91, and 2.44 meters, respectively. The mass of a bale is 0.392 dry Mg. We assume that the sale price of pellets is \$77.16/dry Mg. We assume that unit cost of adjustment of reactor feeding is \$5.51/(dry Mg/minute). We conduct a sensitivity analysis with different values of sale price of pellets and unit cost of adjustment of reactor feeding. We observe that solutions are not sensitive to these parameters.

\subsection{Techno-economic Data}

We use the Biomass Logistics Model (BLM) to estimate input cost data, such as the unit operation cost of pre-processing, handling, and storing used in the model \cite{cafferty2013model}. BLM is an integrated framework that quantifies the impacts of feedstock supply system design on costs, greenhouse gas emissions, and energy consumption. It uses a system dynamics simulation model to simulate the flow of biomass through the whole supply chain, track the change in biomass physical characteristics during operations, and adjust cost/benefit of the engineered system resulting from these operations. BLM estimates the cost of operating a feeding system in \$/hour (see Table \ref{supp-tab:costdata}).

In the case when the capacity of the storage units expands, the corresponding operation cost per minute increases. We use the following scaling rule \cite{tribe1986scale} to estimate the new costs.

\begin{align}
\frac{c_{is}}{c_{isk}} & = \left( \frac{\overline{m}_{i}}{(1+k) \overline{m}_{i}} \right)^{0.6} & \forall i \in \mathcal{I_{M}}, \forall s \in \mathcal{L}, \forall k \in \mathcal{C}
\end{align}

\noindent where 0.6 represents economies of scale. 

\section{Results} \label{sec:results}

In this section, we present the results from implementing the models developed for process design, namely BFFPC and HPC, with and without fractional milling. We compare the performance of BFFPC and HPC by investigating: (\textit{i}) the optimal trajectories; (\textit{ii}) the average reactor feeding rate and its variability; (\textit{iii}) the effect of buffer size expansion and location on system's performance; (\textit{iv}) the cost of operating the system; and ($v$) the minimum time to process all bales.

We assume that a total of 200 bales are to be processed in the feeding system. Of these, 60 are low-moisture bales, 100 are medium-moisture bales, and 40 are high-moisture bales. These bales are processed in HPC following certain {patterns}, which are to process six low-moisture bales, then to process ten medium-moisture bales, and finally, to process four high-moisture bales. This {pattern} is repeated until all 200 bales are processed. This {pattern} is known as 6L,10M,4H X10. BFFPC processes the same number of bales with same moisture distribution. However, bales are not sequenced according to their moisture level. The {feeding pattern} is determined randomly and known as ``random {pattern}'' in this section.

The following is a list of performance metrics we use to compare HPC and BFFPC:

\vspace{-0.5em}
\begin{itemize}
	\itemsep-0.5em
	\item Maximum inventory (\textit{Max. Inv.}) in storage units is $\max_{t} \sum_{s \in \mathcal{L}} \mathbb{M}_{ist}, \forall i \in \mathcal{I_{M}}$.
	\item Average reactor feeding (\textit{Av. Feed.}) is $\frac{\sum_{t \in \mathcal{T}} X_{|\mathcal{I}|t} }{|\mathcal{T}|}$.
	\item Variability of reactor's feeding (\textit{Coef. Var.}) is measured via the coefficient of variation of $X_{|\mathcal{I}|t}$, which is  $\frac{\text{std} \left( X_{|\mathcal{I}|t} \right)}{\left(\frac{\sum_{t \in \mathcal{T}} X_{|\mathcal{I}|t}}{|\mathcal{T}| } \right)}$ where $\text{std}$ is standard deviation operator.
	\item  Cost of operating systems (\textit{Cost}) for HPC and BFFPC is calculated via \eqref{eq:costcalculation} and \eqref{eq:costofbffpc}, respectively.
	\item Minimum time to process 200 bales (\textit{Min. Time}).
\end{itemize}

There are other studies that use the coefficient of variation to measure the variability in production systems, examples of which are given in systems reliability literature \cite{nwanya2016process, li2000production}. Coefficient of variation shows the extent of deviation from average performance in a process. In our case, it indicates how much the feeding rate of a reactor differs from one period to the next.

\subsection{Optimal trajectories of feeding rate of reactor, in-feed rate, and inventory} \label{sec:trajectories}

\begin{figure}[H]
	\begin{subfigure}[b]{\linewidth}
		\centering
		\includegraphics[width=\textwidth]{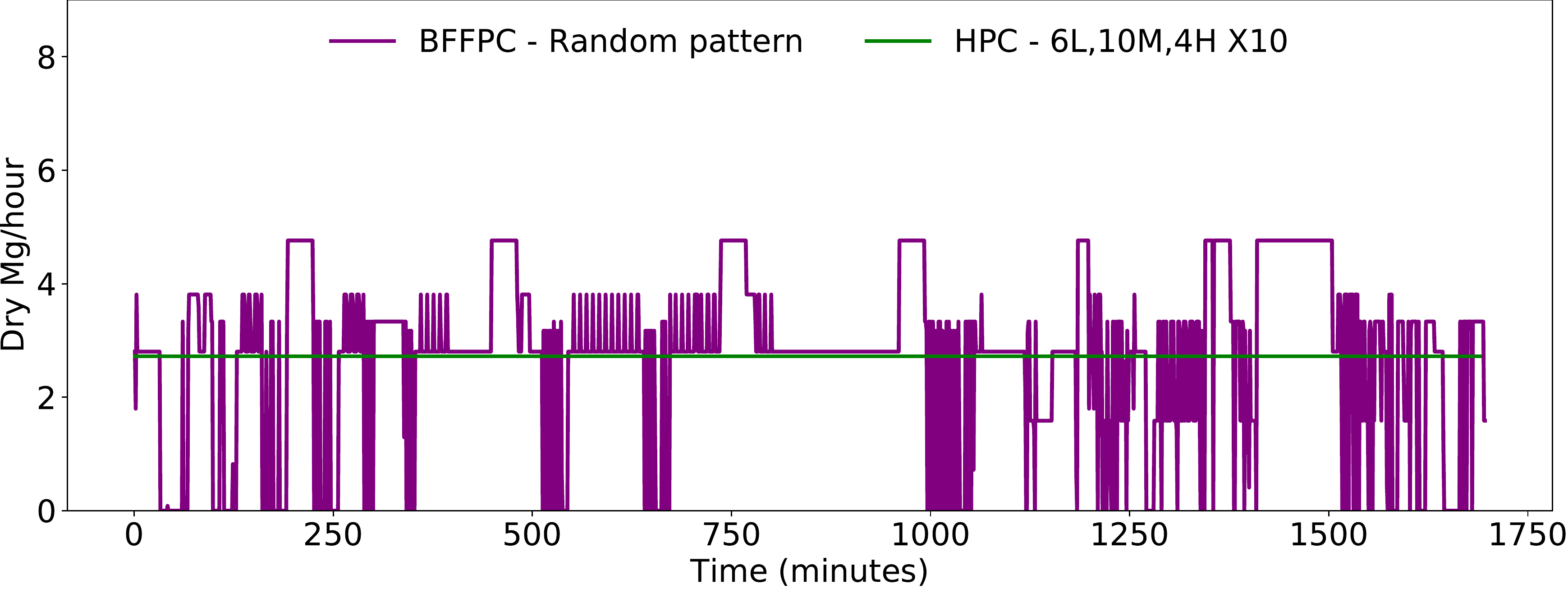}
		\caption{Without fractional milling}
		\label{fig:reactor_trajectory_without}
		\vspace{2ex}
	\end{subfigure}
	\begin{subfigure}[b]{\linewidth}
		\centering
		\includegraphics[width=\textwidth]{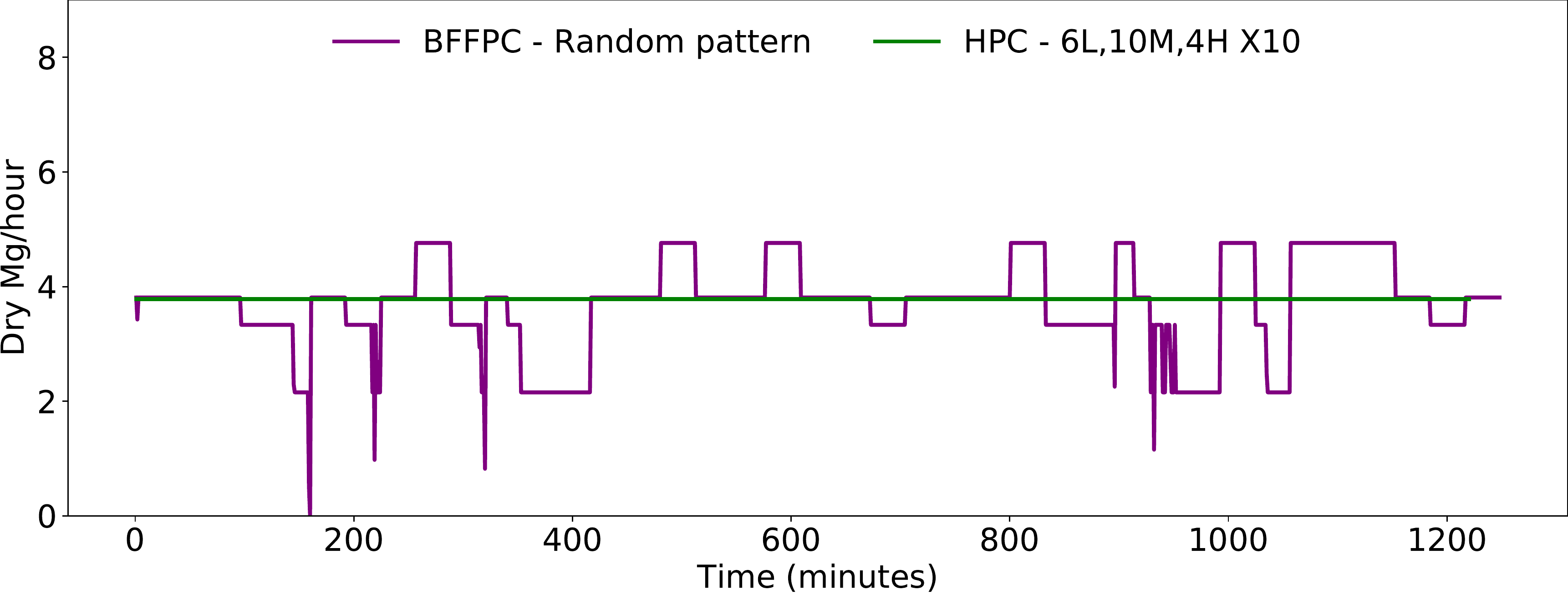}
		\caption{With fractional milling}
		\label{fig:reactor_trajectory_with}
	\end{subfigure}
	\caption{Optimal trajectories of reactor feeding rates for BFFPC (with feeding rate ``random {pattern}'') and HPC (with feeding rate 6L,10M,4H X10)}
	\label{fig:reactor_trajectory}
\end{figure}

In this section, we illustrate the optimal trajectories of a reactor's feeding rate, a system's in-feed rate, and the inventory level. Figure \ref{fig:reactor_trajectory} (a) illustrates the optimal trajectories of reactor's feeding rates for BFFPC and HPC for the process design without fractional milling. Reactor's feeding rate in BFFPC varies from 0 to 5 dry Mg/hour, whereas it is constant in HPC. In the process with fractional milling, BFFPC feeds the reactor at a rate ranging from 2 dry Mg/hour to 5 dry Mg/hour, whereas the HPC feeds the reactor at a constant rate (see Figure \ref{fig:reactor_trajectory} (b)). These results indicate that HPC reduces the variation of a reactor's feeding rate by 100\%, regardless of whether fractional milling is utilized. We further observe that processing times of HPC are slightly shorter than those of BFFPC. It is observed that BFFPC without fractional milling exhibits more variability as compared to BFFPC with fractional milling.

\begin{figure}[H]
	\begin{subfigure}[b]{\linewidth}
		\centering
		\includegraphics[width=\textwidth]{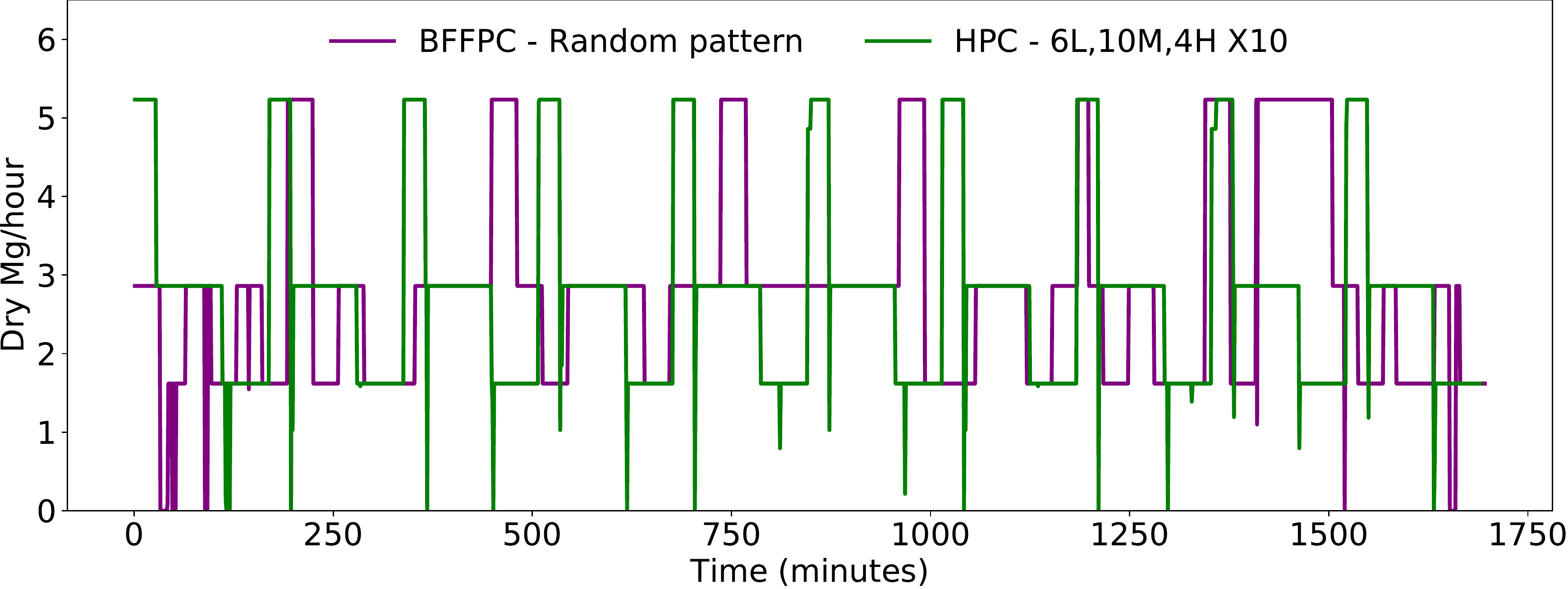}
		\caption{Without fractional milling}
		\label{fig:infeed_trajectory_without}
		\vspace{2ex}
	\end{subfigure}
	\begin{subfigure}[b]{\linewidth}
		\centering
		\includegraphics[width=\textwidth]{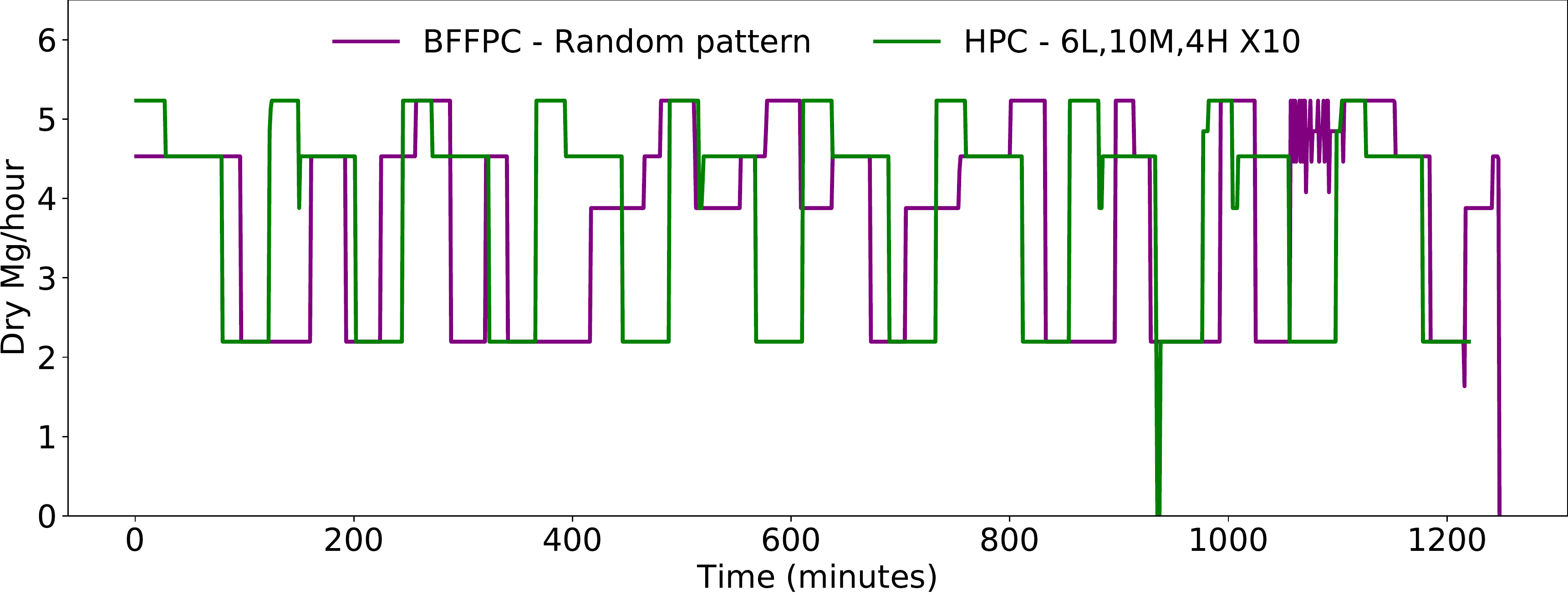}
		\caption{With fractional milling}
		\label{fig:infeed_trajectory_with}
	\end{subfigure}
	\caption{Optimal trajectories of bale in-feed rates for BFFPC (with feeding rate ``random  {pattern}'') and HPC (with feeding rate 6L,10M,4H X10)}
	\label{fig:infeed_trajectory}
\end{figure}

Figure \ref{fig:infeed_trajectory} compares the trajectories of a systems' in-feed rate. We observe that both BFFPC and HPC adjust bale in-feed rates differently. Although the maximum in-feed rate for both pre-processing systems is the same, the sequences of bale fed are different. HPC's in-feed trajectory has a more consistent pattern (same repeated adjustments). If we compare in-feed rate trajectories with and without fractional milling, it is observed that the range of the in-feed rate lies in an interval of 0 dry Mg/hour and 5 dry Mg/hour without fractional milling, this range is  narrower with fractional milling. For instance, we rarely see in-feed rate drops to 0 dry Mg/hour with fractional milling.

\begin{figure}[H]
	\begin{subfigure}[b]{\linewidth}
		\centering
		\includegraphics[width=\textwidth]{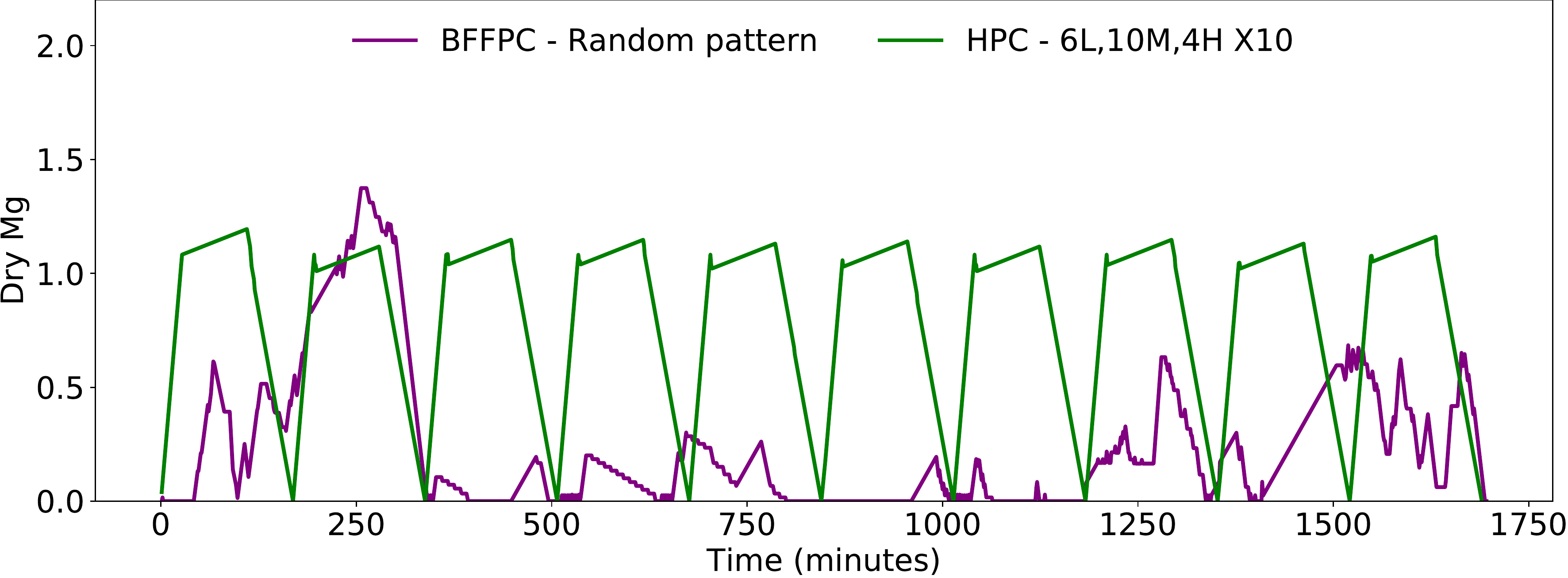}
		\caption{Without fractional milling}
		\label{fig:inventory_bffpc_hpc_without}
		\vspace{2ex}
	\end{subfigure}
	\begin{subfigure}[b]{\linewidth}
		\centering
		\includegraphics[width=\textwidth]{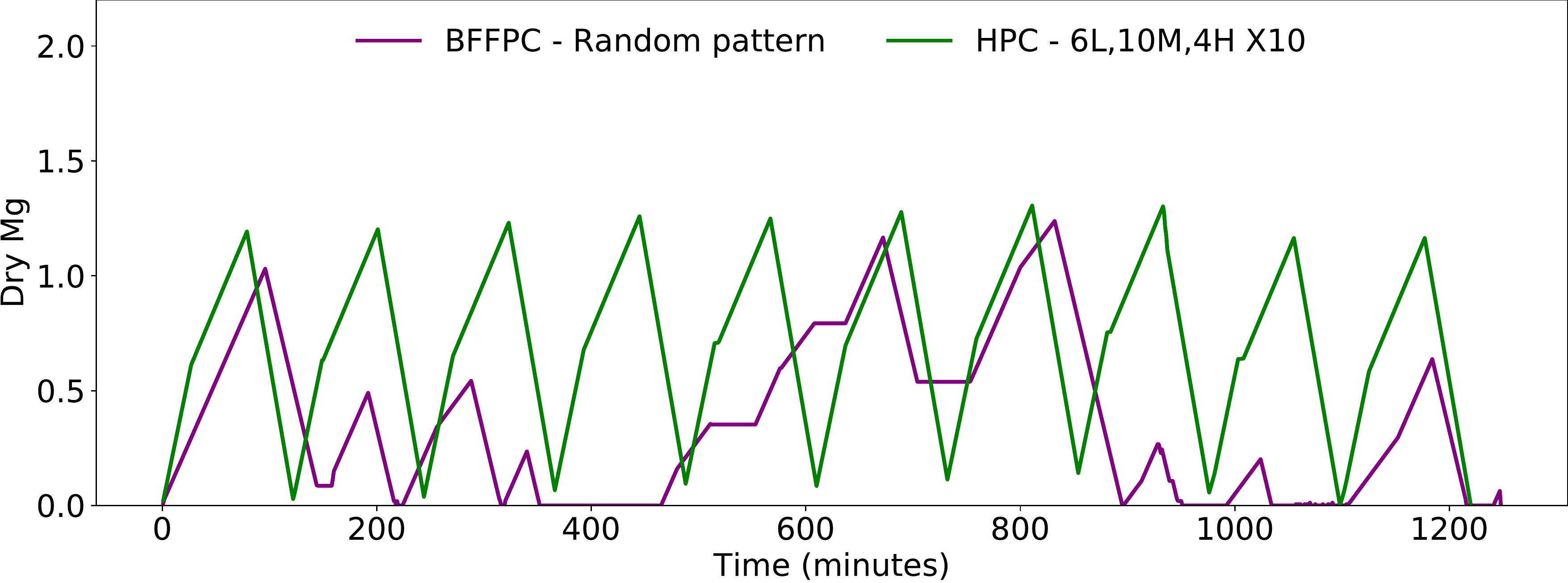}
		\caption{With fractional milling}
		\label{fig:inventory_bffpc_hpc_with}
	\end{subfigure}
	\caption{Optimal trajectories of total inventory in metering and storage bins for BFFPC (with feeding rate ``random {pattern}'') and HPC (with feeding rate 6L,10M,4H X10)}
	\label{fig:inventory_bffpc_hpc} 
\end{figure}

Figure \ref{fig:inventory_bffpc_hpc} reveals optimal trajectories of the total inventory level in metering and storage bins. We observe that the build up of inventory in storage units for BFFPC is irregular. On the other hand, inventory build up for HPC follows a pattern. Biomass accumulates in storage units while processing low- and medium-moisture bales, and biomass is discharged while processing high-moisture bales. The breakpoints observed for HPC correspond to the times when operators switch between bales with different moisture levels. This approach to inventory control minimizes inventory level. The results of Figure \ref{fig:inventory_bffpc_hpc} show that HPC tends to keep more inventory than BFFPC. The higher inventory level leads to a shorter processing time, as well as continuous and smooth feeding of the reactor.

For the {feeding patterns} presented in Figure \ref{fig:inventory_bffpc_hpc}, we calculate the total inventory carried over the entire modeling horizon. It is observed that BFFPC with fractional milling keeps more inventory than BFFPC without fractional milling. However, HPC without fractional milling carries more inventory than HPC with fractional milling. The same conclusions are reached when analyzing HPC and BFFPC with other feeding patterns.

\subsection{Average feeding rate of reactor and variability in feeding}

\begin{figure}[H]
	\begin{subfigure}[b]{0.5\linewidth}
		\centering
		\includegraphics[width=\textwidth]{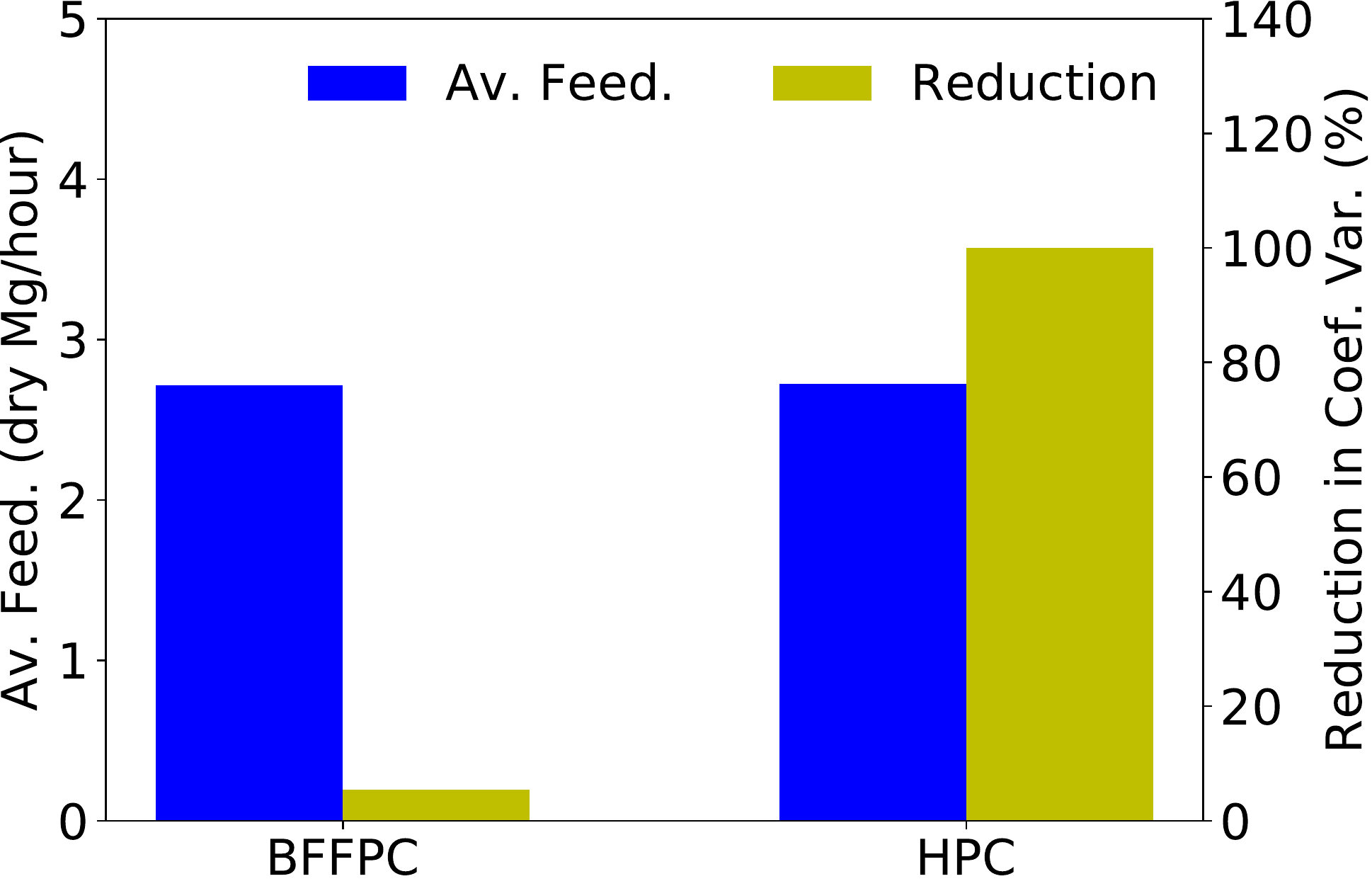}
		\caption{Without fractional milling}
		\label{fig:feedvar_without}
	\end{subfigure}
	\begin{subfigure}[b]{0.5\linewidth}
		\centering
		\includegraphics[width=\textwidth]{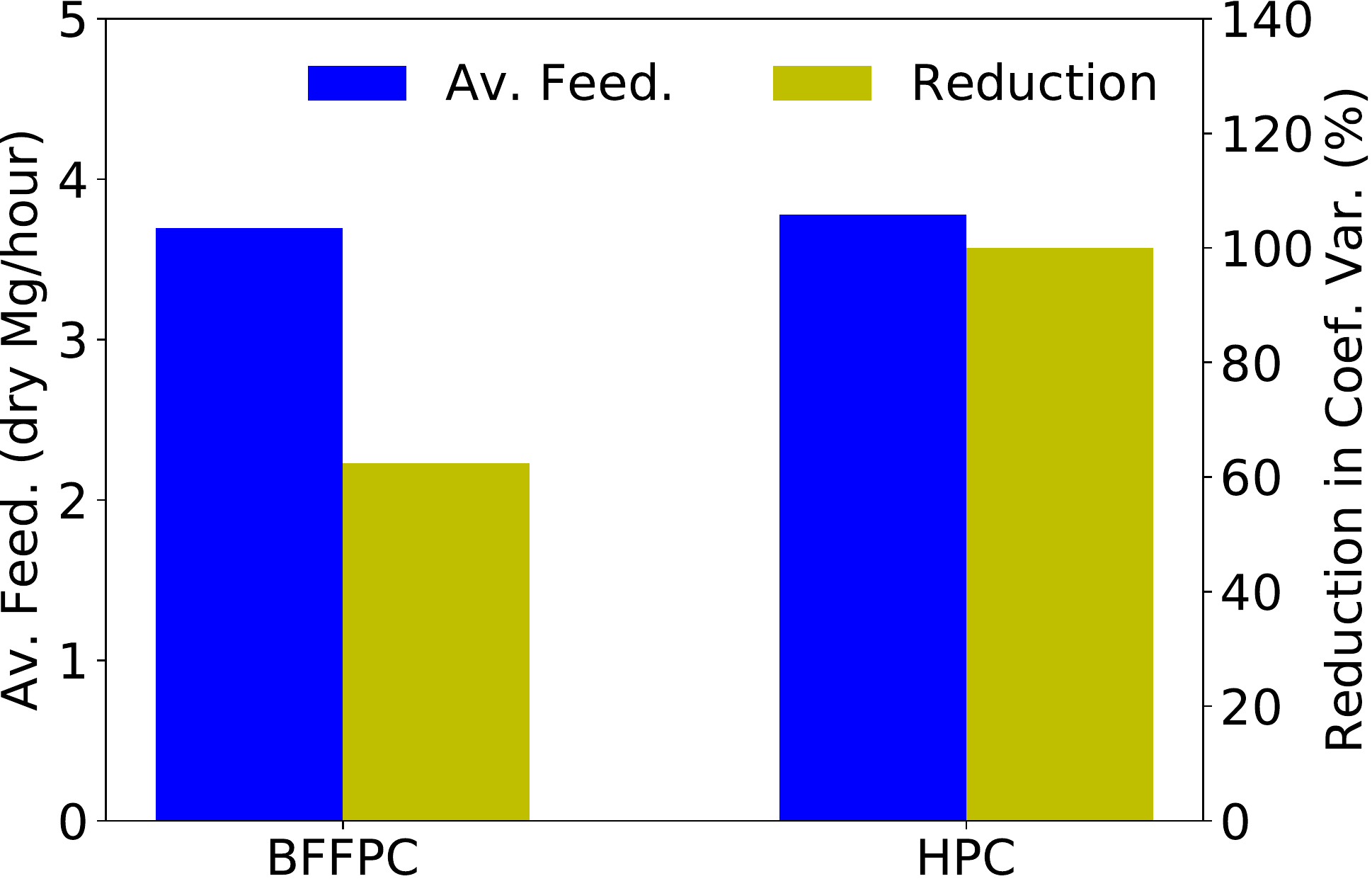}
		\caption{With fractional milling}
		\label{fig:feedvar_with}
	\end{subfigure}
	\caption{Average feeding rate of reactor and reduction in coefficient of variation of reactor feeding rate for BFFPC (with feeding rate ``random {pattern}'') and HPC (with feeding rate 6L,10M,4H X10)}
	\label{fig:feeding_perf}
\end{figure}

In this section, we compare the average feeding rate of the reactor and the variability in feeding. If BFFPC is replaced by HPC in the process design without fractional milling, a reduction in the coefficient of variation  {goes from} 5.48\% {to} 100\%. The average feeding rates {increases from} 2.713 dry Mg/hour {to} 2.723 dry Mg/hour, as shown in Figure \ref{fig:feedvar_without}. If BFFPC is replaced by HPC in the process design with fractional milling, similar improvements are observed in the average feeding rate of reactor and variability in feeding. HPC increases the average feeding rate from 3.696 dry Mg/hour to 3.780 dry Mg/hour, as shown in Figure \ref{fig:feedvar_with}. The reduction in the coefficient of variation in BFFPC and HPC are 62.38\% and 100\%, respectively.

\subsection{Optimal size of buffers to improve the reactor feeding rate}

In this section, we evaluate the impact of buffer size and location to the performance of the feeding system. We remind the reader that the proposed MILP formulation allows capacity expansions of storage units, if economical. We experiment with different feeding patterns in which 60 low-moisture bales, 100 medium-moisture bales, and 40 high-moisture bales are fed to the system in this order (we denote this pattern with 60L,100M,40H). When using the feeding pattern 6L, 10M, 4H X10, the system does not inventory in large amounts; thus, the storage capacity is often not violated. The storage capacity impacts the performance of the system when processing this new pattern.

\begin{figure}[H]
	\begin{subfigure}[b]{0.5\linewidth}
		\centering
		\includegraphics[width=\textwidth]{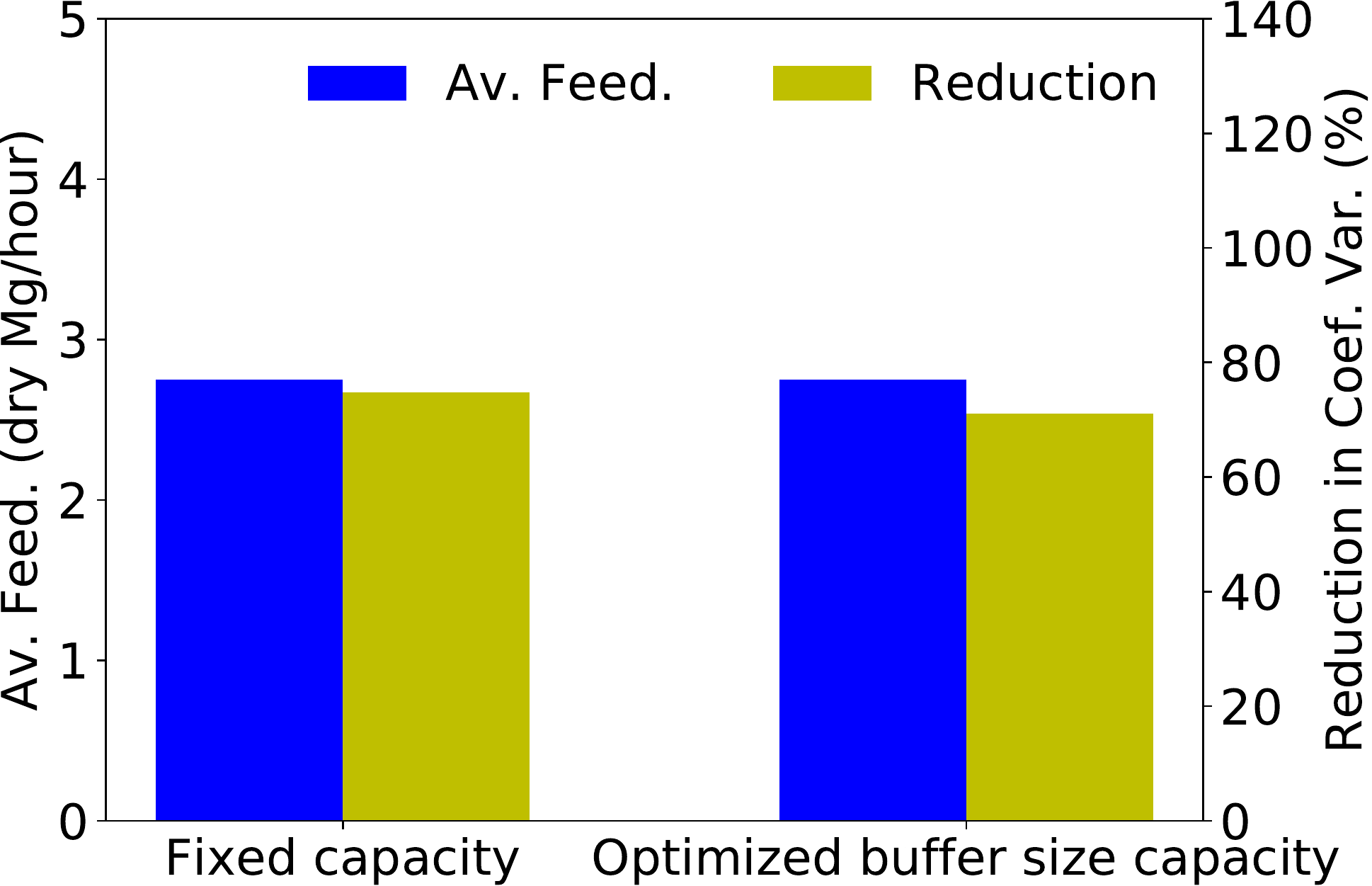}
		\caption{Without fractional milling}
		\label{fig:feedvar_without_buffer}
	\end{subfigure}
	\begin{subfigure}[b]{0.5\linewidth}
		\centering
		\includegraphics[width=\textwidth]{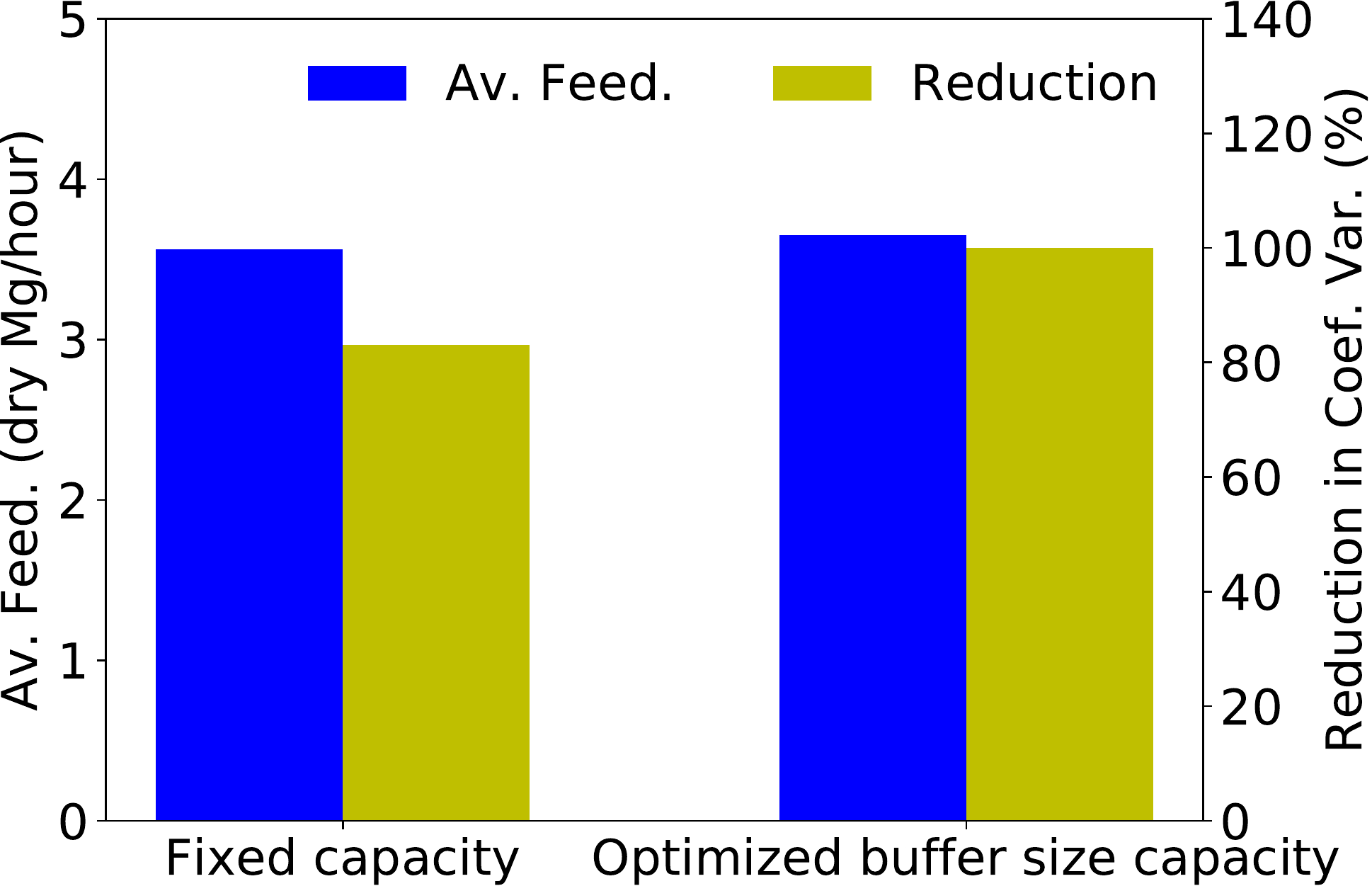}
		\caption{With fractional milling}
		\label{fig:feedvar_with_buffer}
	\end{subfigure}
	\caption{Average feeding rate of reactor and coefficient of variation in reactor feeding rate with fixed and optimized capacities of storage units}
	\label{fig:feedvar_buffer}
\end{figure}

\begin{figure}[H]
	\begin{subfigure}[b]{0.5\linewidth}
		\centering
		\includegraphics[width=\textwidth]{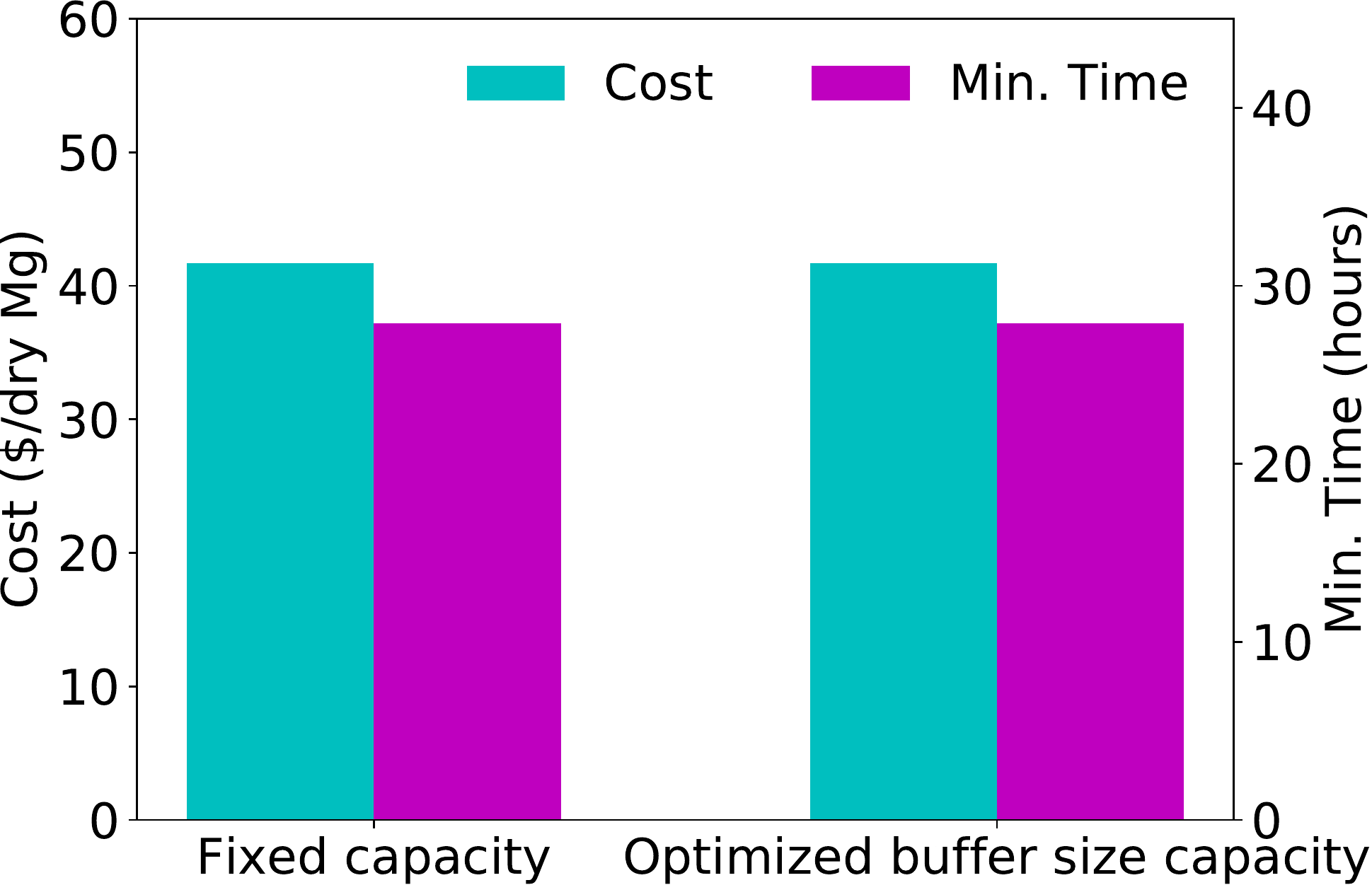}
		\caption{Without fractional milling}
		\label{fig:cost_time_without_capanalysis}
	\end{subfigure}
	\begin{subfigure}[b]{0.5\linewidth}
		\centering
		\includegraphics[width=\textwidth]{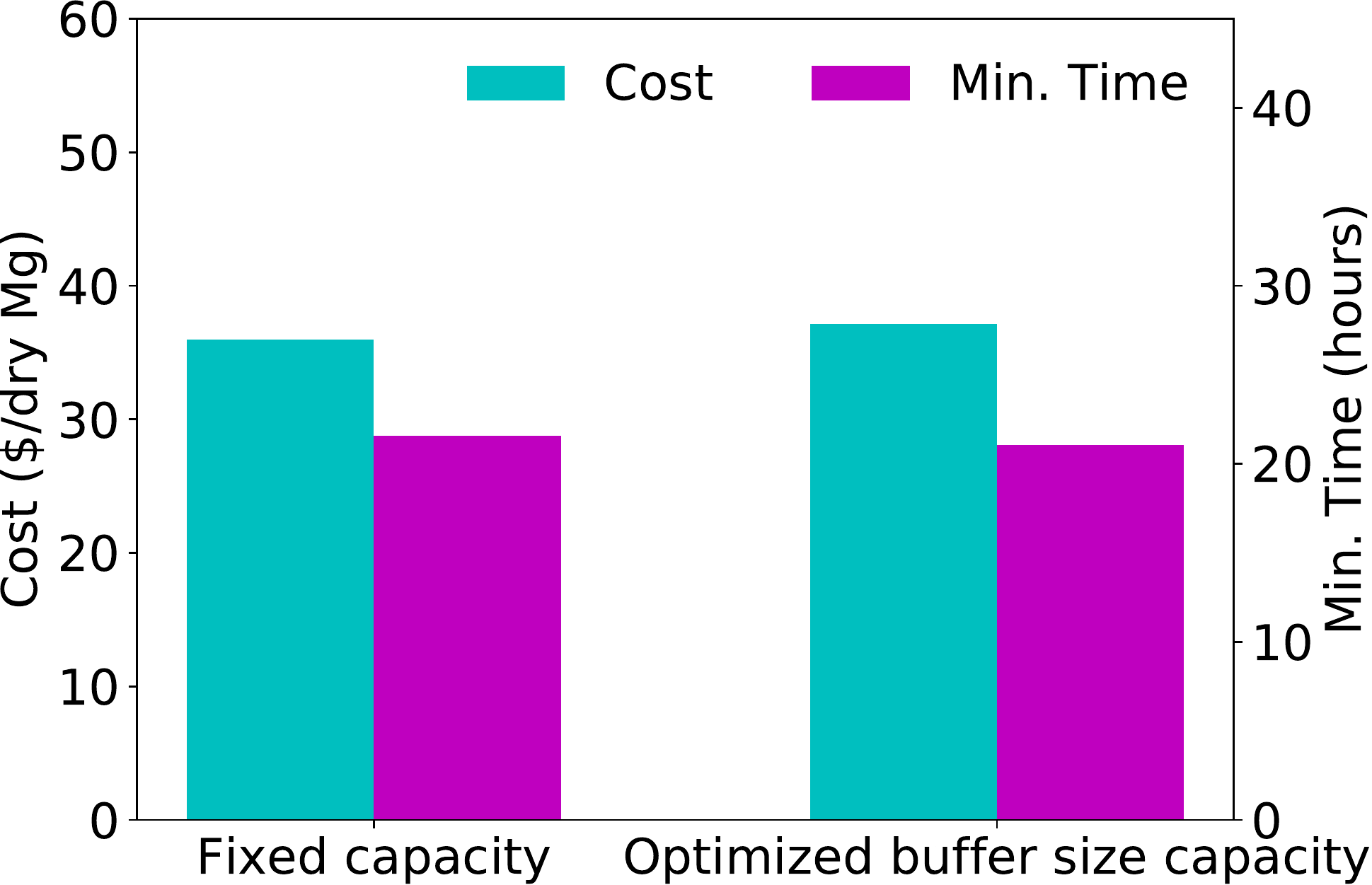}
		\caption{With fractional milling}
		\label{fig:cost_time_with_capanalysis}
	\end{subfigure}
	\caption{Cost of operating system and minimum time to process all bales with fixed and optimized capacities of storage units}
	\label{fig:cost_time_capanalysis}
\end{figure}

In Figures \ref{fig:feedvar_buffer} and  \ref{fig:cost_time_capanalysis}, ``optimized buffer size capacity" refers to the model that allows expansion (as needed) of the capacity of storage units, and ``fixed capacity" refers to the model that does not allow the expansion of the capacity of storage units. We note that the system without fractional milling does not need expansion in storage units; and thus optimized buffer size capacity equals original capacities. In an optimal solution of the system with fractional milling, buffer size capacities are expanded by 100\%.

The results of Figures \ref{fig:feedvar_with_buffer} and Figure \ref{fig:cost_time_with_capanalysis} show that the HPC system improves its performance when processing the feeding pattern 60L, 100M, 40H and additional  storage capacity is available. In this case, the process throughput is 3.649 dry Mg/hour, the reduction in the coefficient of variation is 100\%, the average unit cost is \$37.15/dry Mg, and the processing time is 21.07 hours. If the existing storage capacity (of 4.54 dry Mg and 24.07 cubic meter) in the process design did not increase, the average throughput would be lower (3.564 dry Mg/hour), the reduction in coefficient of variation would be lower (83.06\%), the average unit cost would be lower (\$35.97/dry Mg), and the process time would be higher (21.57 hours). The additional investment cost is the main reason why the average unit cost is higher in a system that uses additional storage capacity than a system that does not. These additional costs outweigh the reductions of processing costs due to shorter processing times. The impact of buffer size is minimal to the process design without fractional milling. This is because system feeding is relatively slower and it does not require extra storage capacity.

\subsection{Cost of operating the system}

Table \ref{tab:costtime} summarizes the estimated costs of operating the system and the minimum time required to process 200 bales for BFFPC and HPC. We observe that HPC reduces costs by 0.36\% and processing time by 0.35\% over BFFPC in a process without fractional milling. When fractional milling is utilized, HPC's economical advantages over BFFPC increase. HPC saves 2.22\% in costs and 2.24\% in processing time. These results also indicate improvements of processing time and costs due to fractional milling.

\begin{table}[H]
	\footnotesize
	\centering
	\begin{tabular}{llcc}
		\toprule
		\textbf{Structure} & \textbf{Control} & \textbf{Cost (\$/dry Mg)} & \textbf{Min. Time (hours)} \\
		\midrule
		\multirow{2}{*}{Without fractional milling} & BFFPC & 42.32 & 28.27 \\
		& HPC & 42.16 & 28.17 \\
		\midrule
		\multirow{2}{*}{With fractional milling} & BFFPC & 34.71 & 20.80 \\
		& HPC & 33.94 & 20.33 \\
		\bottomrule
	\end{tabular}
	\caption{Cost of operating the system and minimum time to process 200 bales for BFFPC (with feeding rate ``random {pattern}'') and HPC (with feeding rate 6L,10M,4H X10)}
	\label{tab:costtime}
\end{table}

\section{Conclusion and Future Research Directions} \label{sec:conclusion}

The variations in feedstock characteristics, such as feedstock moisture and particle size distribution lead to an inconsistent flow of feedstock from the biomass pre-processing system to a rector in-feed system. In this research, we developed an optimal process control method for a biomass pre-processing system comprised of milling and densification operations to provide consistent flow of feedstock to a reactor's throat. Our case study based on switchgrass finds that the proposed optimal process control reduces the variation of reactor feeding rate by 100\% without increasing the cost of a biomass pre-processing system operating if the biomass moisture ranges from 10 to 25\%. This variability reduction is possible via identification of optimal bale sequencing, equipment in-feed rate in the process, and process buffer location and size in the biomass pre-processing system. A biorefinery can adapt proposed control to achieve its design capacity.

Future research will focus on integrating the outcomes of Discrete Element Method (DEM) models with optimization models. DEM models are typically used to simulate the mechanical and flow behavior of feedstocks in equipment for given feedstock characteristics and equipment design. The results from these simulations can be used to develop functional relationships between feedstock characteristics (e.g., type, moisture content, particle size distribution, density), equipment setting, and biomass flow. These functional relationships can be incorporated into the optimization model to control the flow of biomass to the reactor. For example, biomass type, moisture level and mill speed in grinders impact biomass density. The density of biomass impacts its volume, and as a consequence, the storage capacity of the metering bin. Thus, outcomes of DEM models can be used to develop functional relationships between biomass type, biomass moisture level, mill speed in grinders and density to be included in optimization. There are a few DEM models for biomass feedstocks, such as, Guo et al. \cite{guo2020discrete} which focuses on switchgrass and Xia et al. \cite{xia2019discrete} which focuses on pinewood. We will incorporate the outcomes of similar models to optimization to improve the accuracy of our results.
\section{Credit author statement}
Fikri Kucuksayacigil: Develop and implementation of model utilizing cplex/Gurobi, prepare results, writing original draft. Mohammad Roni: Conceptualization, develop methodology, scoping the study, data collection, writing: prepare final manuscript Sandra D. Eksioglu: Develop methodology, supervision, writing- review and editing.  Qiushi Chen:Writing - review \& editing. 
\section{Declaration of competing interest}
The authors declare that they have no known competing financial interests or personal relationships that could have appeared to influence the work reported in this paper.
\section{Acknowledgment}

The funding of this work was provided by the U.S. Department of Energy, Office of Energy Efficiency and Renewable Energy, Bioenergy Technologies Office under award Number DE- EE0008255 and Department of Energy Idaho Operations Office Contract No. DE-AC07-05ID14517. The authors greatly thank INL staff in the DOE Biomass Feedstock National User Facility for providing technical data for this study. Specifically, we acknowledge Neal A. Yancey and Jaya S. Tumuluru, who are the individual INL contributors to this research who provided helpful comments, data, and other forms of support for this analysis. The views expressed herein are those of the author only, and do not necessarily represent the views of DOE or the U.S. Government.

%
%

\bibliographystyle{elsarticle-num}
\bibliography{References_paper}

\begin{thebibliography}{10}
\expandafter\ifx\csname url\endcsname\relax
  \def\url#1{\texttt{#1}}\fi
\expandafter\ifx\csname urlprefix\endcsname\relax\def\urlprefix{URL }\fi
\expandafter\ifx\csname href\endcsname\relax
  \def\href#1#2{#2} \def\path#1{#1}\fi

\bibitem{us2016biorefinery}
United States Department of Energy, Chicago, Illinois, Biorefinery Optimization
  Workshop Summary Report (2016).

\bibitem{bell2005challenges}
T.~A. Bell, Challenges in the scale-up of particulate processes—{A}n
  industrial perspective, Powder Technology 150~(2) (2005) 60--71.

\bibitem{merrow2000developing}
E.~W. Merrow, Developing technology-problems and progress in particle
  processing-careful planning and design forestall many of the problems
  inherent in commercializing processes that use solids processing, Chemical
  Innovation 30~(1) (2000) 34--41.

\bibitem{williams2016sources}
C.~L. Williams, T.~L. Westover, R.~M. Emerson, J.~S. Tumuluru, C.~Li, Sources
  of biomass feedstock variability and the potential impact on biofuels
  production, BioEnergy Research 9~(1) (2016) 1--14.

\bibitem{westover2018biomass}
T.~Westover, D.~S. Hartley, Biomass handling and feeding, Advances in Biofuels
  and Bioenergy (2018) 86.

\bibitem{hirtzer2017}
M.~Hirtzer,
  \href{https://www.reuters.com/article/us-dowdupontethanol/dupont-to-sell-
  \newline cellulosic-ethanol-plant-in-blowto-biofuel-idUSKBN1D22T5}{DuPont to
  sell cellulosic ethanol plant in blow to biofuel} (2017).
\newline\urlprefix\url{https://www.reuters.com/article/us-dowdupontethanol/dupont-to-sell-
  \newline cellulosic-ethanol-plant-in-blowto-biofuel-idUSKBN1D22T5}

\bibitem{mamun2020supply}
S.~Mamun, J.~K. Hansen, M.~S. Roni, Supply, operational, and market risk
  reduction opportunities: Managing risk at a cellulosic biorefinery, Renewable
  and Sustainable Energy Reviews 121 (2020) 109677.

\bibitem{hess2020advancements}
J.~R. Hess, A.~E. Ray, T.~G. Rials, Advancements in Biomass Feedstock
  Preprocessing: Conversion Ready Feedstocks, Frontiers Media SA, 2020.

\bibitem{dooley2013woody}
J.~H. Dooley, D.~N. Lanning, C.~J. Lanning, Woody biomass size reduction with
  selective material orientation, Biofuels 4~(1) (2013) 35--43.

\bibitem{ketterhagen2009predicting}
W.~R. Ketterhagen, J.~S. Curtis, C.~R. Wassgren, B.~C. Hancock, Predicting the
  flow mode from hoppers using the discrete element method, Powder Technology
  195~(1) (2009) 1--10.

\bibitem{crawford2016effects}
N.~C. Crawford, N.~Nagle, D.~A. Sievers, J.~J. Stickel, The effects of physical
  and chemical preprocessing on the flowability of corn stover, Biomass and
  Bioenergy 85 (2016) 126--134.

\bibitem{wu2011physical}
M.~Wu, D.~Schott, G.~Lodewijks, Physical properties of solid biomass, Biomass
  and Bioenergy 35~(5) (2011) 2093--2105.

\bibitem{marino2017data}
D.~Marino, K.~Amarasinghe, M.~Anderson, N.~Yancey, Q.~Nguyen, K.~Kenney,
  M.~Manic, Data driven decision support for reliable biomass feedstock
  preprocessing, in: 2017 Resilience Week (RWS), IEEE, 2017, pp. 97--102.

\bibitem{marino2018interpretable}
D.~L. Marino, M.~Anderson, K.~Kenney, M.~Manic, Interpretable data-driven
  modeling in biomass preprocessing, in: 2018 11th International Conference on
  Human System Interaction (HSI), IEEE, 2018, pp. 291--297.

\bibitem{jaramillo2018mass}
I.~Jaramillo, A.~Sanchez, Mass flow dynamic modeling and residence time control
  of a continuous tubular reactor for biomass pretreatment, ACS Sustainable
  Chemistry \& Engineering 6~(7) (2018) 8570--8577.

\bibitem{numbi2015systems}
B.~P. Numbi, X.~Xia, Systems optimization model for energy management of a
  parallel {HPGR} crushing process, Applied Energy 149 (2015) 133--147.

\bibitem{hartley2020effect}
D.~S. Hartley, D.~N. Thompson, L.~M. Griffel, Q.~A. Nguyen, M.~S. Roni, Effect
  of biomass properties and system configuration on the operating effectiveness
  of biomass to biofuel systems, ACS Sustainable Chemistry \& Engineering
  8~(19) (2020) 7267--7277.

\bibitem{das2011volume}
D.~Das, A.~Roy, S.~Kar, A volume flexible economic production lot-sizing
  problem with imperfect quality and random machine failure in fuzzy-stochastic
  environment, Computers \& Mathematics with Applications 61~(9) (2011)
  2388--2400.

\bibitem{li2000production}
J.~Li, S.~M. Meerkov, Production variability in manufacturing systems:
  Bernoulli reliability case, Annals of Operations Research 93~(1-4) (2000)
  299--324.

\bibitem{assaf2014analytical}
R.~Assaf, M.~Colledani, A.~Matta, Analytical evaluation of the output
  variability in production systems with general markovian structure, OR
  Spectrum 36~(3) (2014) 799--835.

\bibitem{colledani2008analysis}
M.~Colledani, A.~Matta, T.~Tolio, Analysis of the production variability in
  manufacturing lines, in: Engineering Systems Design and Analysis, Vol. 48357,
  2008, pp. 381--390.

\bibitem{tempelmeier2003practical}
H.~Tempelmeier, Practical considerations in the optimization of flow production
  systems, International Journal of Production Research 41~(1) (2003) 149--170.

\bibitem{tan2009analysis}
B.~Tan, S.~B. Gershwin, Analysis of a general markovian two-stage
  continuous-flow production system with a finite buffer, International Journal
  of Production Economics 120~(2) (2009) 327--339.

\bibitem{turki2013perturbation}
S.~Turki, S.~Hennequin, N.~Sauer, Perturbation analysis for continuous and
  discrete flow models: a study of the delivery time impact on the optimal
  buffer level, International Journal of Production Research 51~(13) (2013)
  4011--4044.

\bibitem{hosseini2017simulation}
B.~Hosseini, B.~Tan, Simulation and optimization of continuous-flow production
  systems with a finite buffer by using mathematical programming, IISE
  Transactions 49~(3) (2017) 255--267.

\bibitem{under2014feedstock}
J.~J. Jacobson, M.~S. Roni, K.~G. Cafferty, K.~Kenney, E.~Searcy, J.~Hansen,
  Feedstock Supply System Design and Analysis ``The Feedstock Logistics Design
  Case for Multiple Conversion Pathways", Idaho National Laboratory, Idaho, USA
  (2014).

\bibitem{kenney2013feedstock}
K.~Kenney, K.~G. Cafferty, J.~J. Jacobson, I.~J. Bonner, G.~L. Gresham, W.~A.
  Smith, D.~N. Thompson, V.~S. Thompson, J.~S. Tumuluru, N.~Yancey, Feedstock
  supply system design and economics for conversion of lignocellulosic biomass
  to hydrocarbon fuels: Conversion pathway: Biological conversion of sugars to
  hydrocarbons the 2017 design case, Tech. rep., Idaho National Laboratory
  (INL) (2013).

\bibitem{kyriakides2020dynamic}
A.-S. Kyriakides, A.~I. Papadopoulos, P.~Seferlis, I.~Hassan, Dynamic modelling
  and control of single, double and triple effect absorption refrigeration
  cycles, Energy 210 (2020) 118529.

\bibitem{arnaudo2020techno}
M.~Arnaudo, M.~Topel, B.~Laumert, Techno-economic analysis of demand side
  flexibility to enable the integration of distributed heat pumps within a
  swedish neighborhood, Energy 195 (2020) 117012.

\bibitem{bezanson2017julia}
J.~Bezanson, A.~Edelman, S.~Karpinski, V.~B. Shah,
  \href{https://doi.org/10.1137/141000671}{Julia: A fresh approach to numerical
  computing}, SIAM Review 59~(1) (2017) 65--98.
\newline\urlprefix\url{https://doi.org/10.1137/141000671}

\bibitem{DunningHuchetteLubin2017}
I.~Dunning, J.~Huchette, M.~Lubin, Jump: A modeling language for mathematical
  optimization, SIAM Review 59~(2) (2017) 295--320.
\newblock \href {https://doi.org/10.1137/15M1020575}
  {\path{doi:10.1137/15M1020575}}.

\bibitem{gurobi2020}
L.~Gurobi~Optimization, \href{http://www.gurobi.com}{Gurobi optimizer reference
  manual} (2020).
\newline\urlprefix\url{http://www.gurobi.com}

\bibitem{mani2004grinding}
S.~Mani, L.~G. Tabil, S.~Sokhansanj, Grinding performance and physical
  properties of wheat and barley straws, corn stover and switchgrass, Biomass
  and Bioenergy 27~(4) (2004) 339--352.

\bibitem{tumuluru2019biomass}
J.~S. Tumuluru, D.~J. Heikkila, Biomass grinding process optimization using
  response surface methodology and a hybrid genetic algorithm, Bioengineering
  6~(1) (2019) 12.

\bibitem{tumuluru2017}
J.~S. Tumuluru,
  \href{https://www.energy.gov/sites/prod/files/2017/05/f34/fsl-tumuluru-1222.pdf}{Biomass
  Engineering: Size reduction, drying and densification of high moisture
  biomass} (2017).
\newline\urlprefix\url{https://www.energy.gov/sites/prod/files/2017/05/f34/fsl-tumuluru-1222.pdf}

\bibitem{yancey2013drying}
N.~A. Yancey, J.~S. Tumuluru, C.~T. Wright, Drying, grinding and pelletization
  studies on raw and formulated biomass feedstock's for bioenergy applications,
  Journal of Biobased Materials and Bioenergy 7~(5) (2013) 549--558.

\bibitem{cafferty2013model}
K.~G. Cafferty, D.~J. Muth, J.~J. Jacobson, K.~M. Bryden, Model based biomass
  system design of feedstock supply systems for bioenergy production, in: ASME
  2013 International Design Engineering Technical Conferences and Computers and
  Information in Engineering Conference, American Society of Mechanical
  Engineers Digital Collection, 2013.

\bibitem{tribe1986scale}
M.~Tribe, R.~Alpine, Scale economies and the “0.6 rule”, Engineering Costs
  and Production Economics 10~(1) (1986) 271--278.

\bibitem{nwanya2016process}
S.~Nwanya, C.~Achebe, O.~Ajayi, C.~Mgbemene, Process variability analysis in
  make-to-order production systems, Cogent Engineering 3~(1) (2016) 1269382.

\bibitem{guo2020discrete}
Y.~Guo, Q.~Chen, Y.~Xia, T.~Westover, S.~Eksioglu, M.~Roni, Discrete element
  modeling of switchgrass particles under compression and rotational shear,
  Biomass and Bioenergy 141 (06 2020).

\bibitem{xia2019discrete}
Y.~Xia, Z.~Lai, T.~Westover, J.~Klinger, H.~Huang, Q.~Chen, Discrete element
  modeling of deformable pinewood chips in cyclic loading test, Powder
  Technology 345 (2019) 1--14.

\end{thebibliography}

\end{document}


\setlength{\abovedisplayskip}{-12pt}
\setlength{\belowdisplayskip}{8pt}
\setlength{\abovedisplayshortskip}{-12pt}
\setlength{\belowdisplayshortskip}{8pt}

\renewcommand{\thefigure}{S\arabic{figure}}
\renewcommand{\thetable}{S\arabic{table}}

\begin{frontmatter}

\title{Optimal Control of Feedstock Pre-processing to Handle Variations in Feedstock Characteristics and Reactor In-Feed Rate}

\author[a]{Fikri Kucuksayacigil}
\author[a]{Mohammad Roni}
\author[b]{Sandra D. Eksioglu}
\author[c]{Qiushi Chen}
\address[a]{Energy and Environment Scienve \& Technology Department, Idaho National Laboratory, Idaho Falls, ID}
\address[b]{Industrial Engineering Department, University of Arkansas, Fayetteville, AR}
\address[c]{Glenn Department of Civil Engineering, Clemson University, Clemson, SC}

\author{}

\address{}

\end{frontmatter}

\newpage

\section{Visual Description of Equipment in the Process}

\subsection{PDU Instrumentation Diagram}



\begin{center}
	\begin{sideways}
		\begin{minipage}{1.3\linewidth}
			\includegraphics[width=\linewidth]{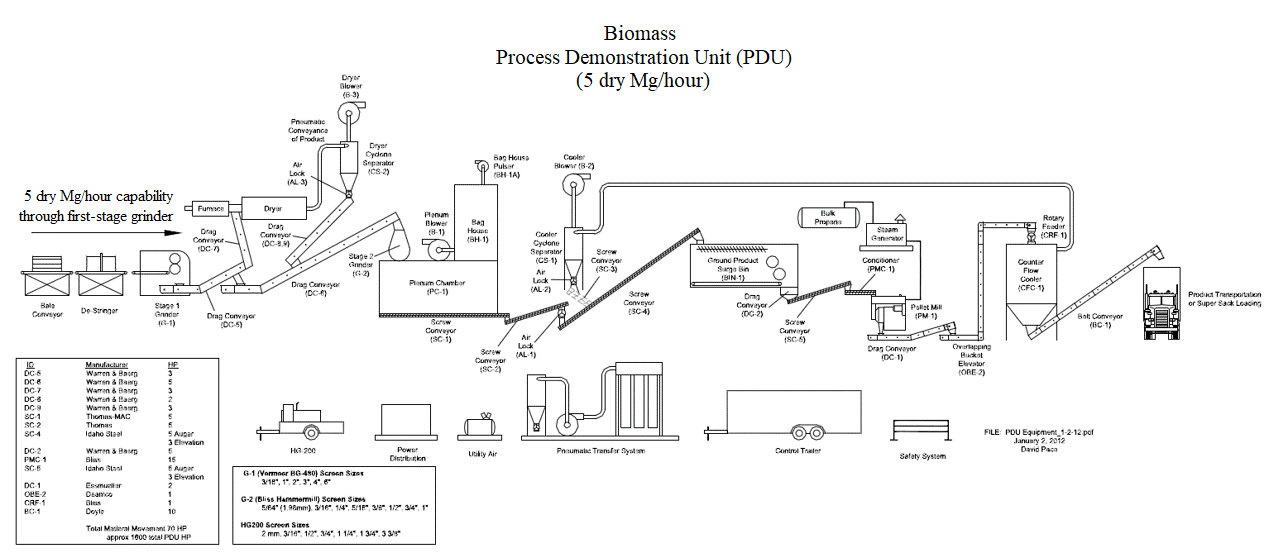}
			\captionof{figure}{PDU technical diagram}
			\vspace{0.2cm}
			\label{fig:pdutechnical}
		\end{minipage}
	\end{sideways}
\end{center}

\subsection{First-stage Grinder}

\begin{figure}[H]
	\centering
	\includegraphics[width=0.9\linewidth]{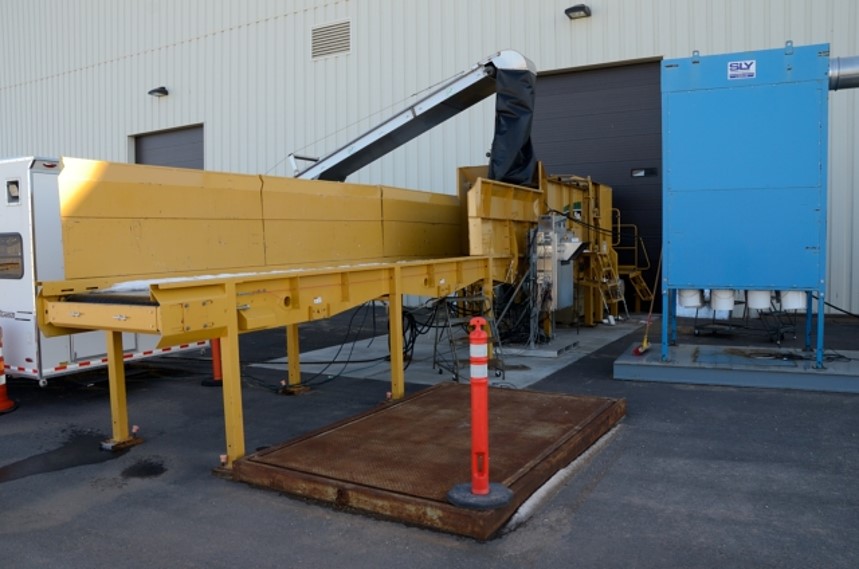}
	\caption{First-stage grinder}
	\label{fig:grinder1}
\end{figure}

\begin{itemize}
	\itemsep-0.5em
	\item This is a Vermeer Prototype Horizontal Bale Grinder.
	\item It is powered by two electric motors (each with 149.14 kW).
	\item It has 2 rotating drums with 96 swinging hammers on each drum.
	\item It can process up to 1.22 m x 1.22 m x 2.44 m rectangular bales.
	\item It can also process loose or baled woody material, but it can not process whole logs.
	\item Its design capacity for dry corn stover is 5 dry Mg/hour.
	\item Various screen sizes are available for this grinder such as 4.76 mm, 25.4 mm, 50.8 mm, 76.2 mm, 101.6 mm, and 152.4 mm.
\end{itemize}

\subsection{Second-stage Grinder}

\begin{figure}[H]
	\centering
	\includegraphics[width=0.9\linewidth]{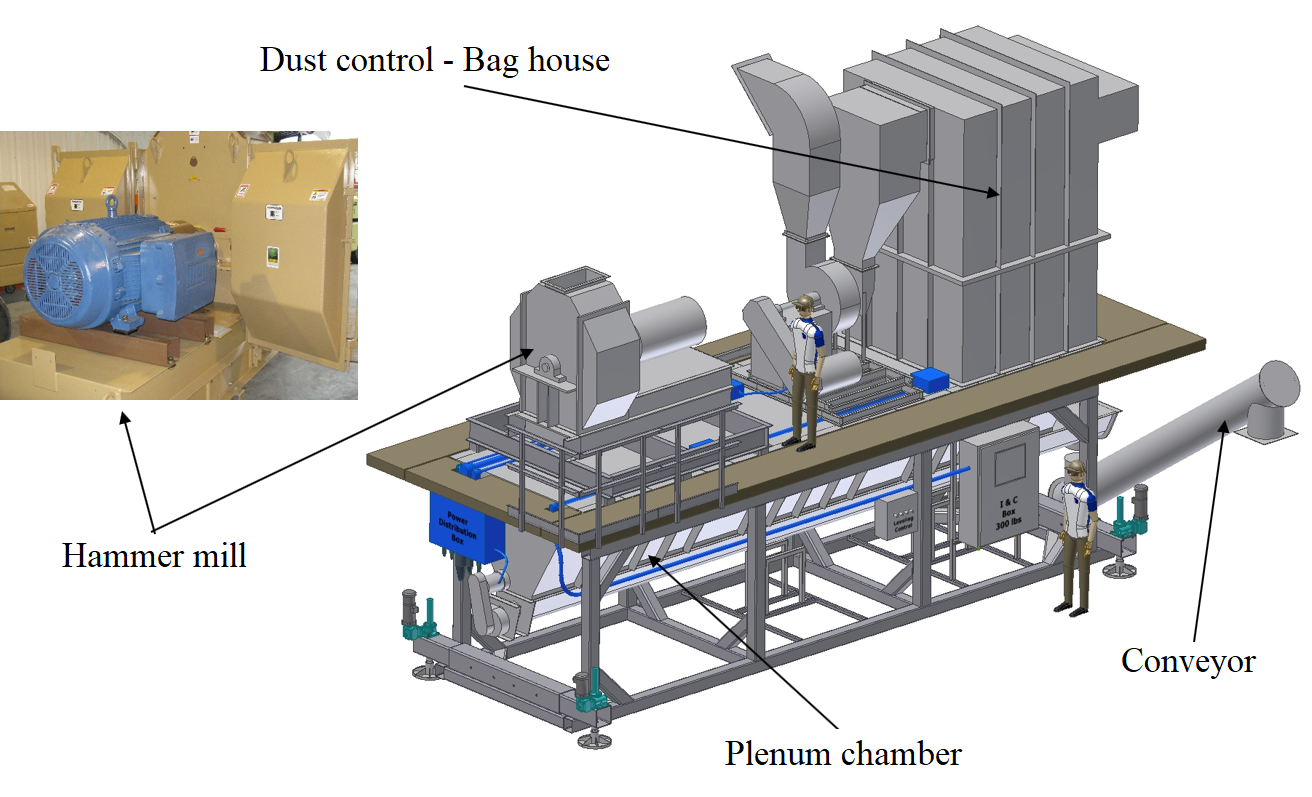}
	\caption{Second-stage grinder}
	\label{fig:grinder2-1}
\end{figure}

\begin{figure}[H]
	\centering
	\includegraphics[width=0.9\linewidth]{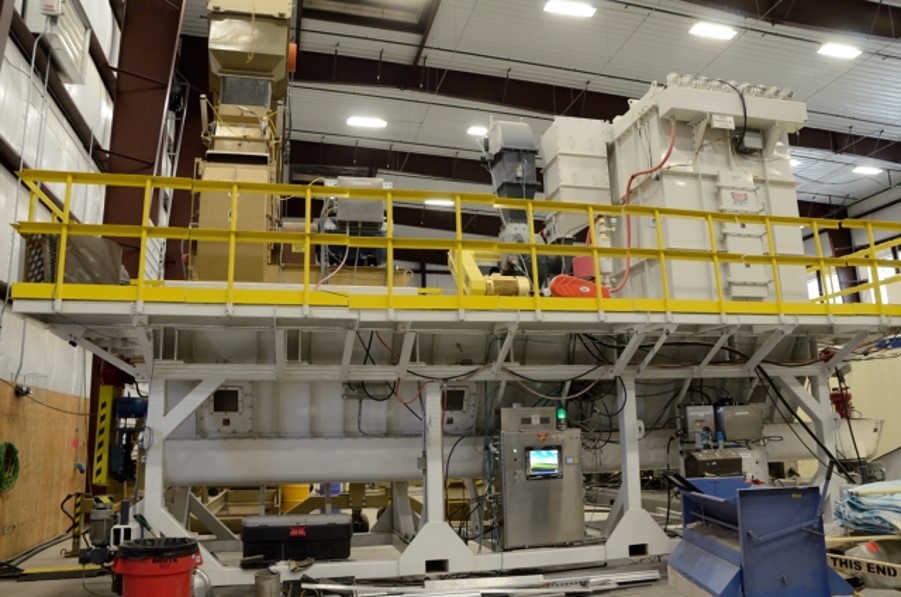}
	\caption{Real picture of second-stage grinder }
	\label{fig:grinder2}
\end{figure}

\begin{itemize}
	\itemsep-0.5em
	\item This is a Bliss Eliminator Series Hammer mill.
	\item It is powered by a 11.86 kW electric motor.
	\item It has 40 swinging hammers.
	\item It has a design capacity for dry corn stover of 5 dry Mg/hour.
	\item Various screen sizes are available for this grinder at the INL PDU, such as 1.98 mm, 4.76 mm, 6.35 mm, 7.94 mm, 9.52 mm, 11.11 mm, 12.7 mm, 15.87 mm, 19.05 mm, 25.4 mm, and 50.8 mm. Other screens are available from Bliss.
\end{itemize}

\subsection{Densification}

\begin{figure}[H]
	\centering
	\includegraphics[width=0.9\linewidth]{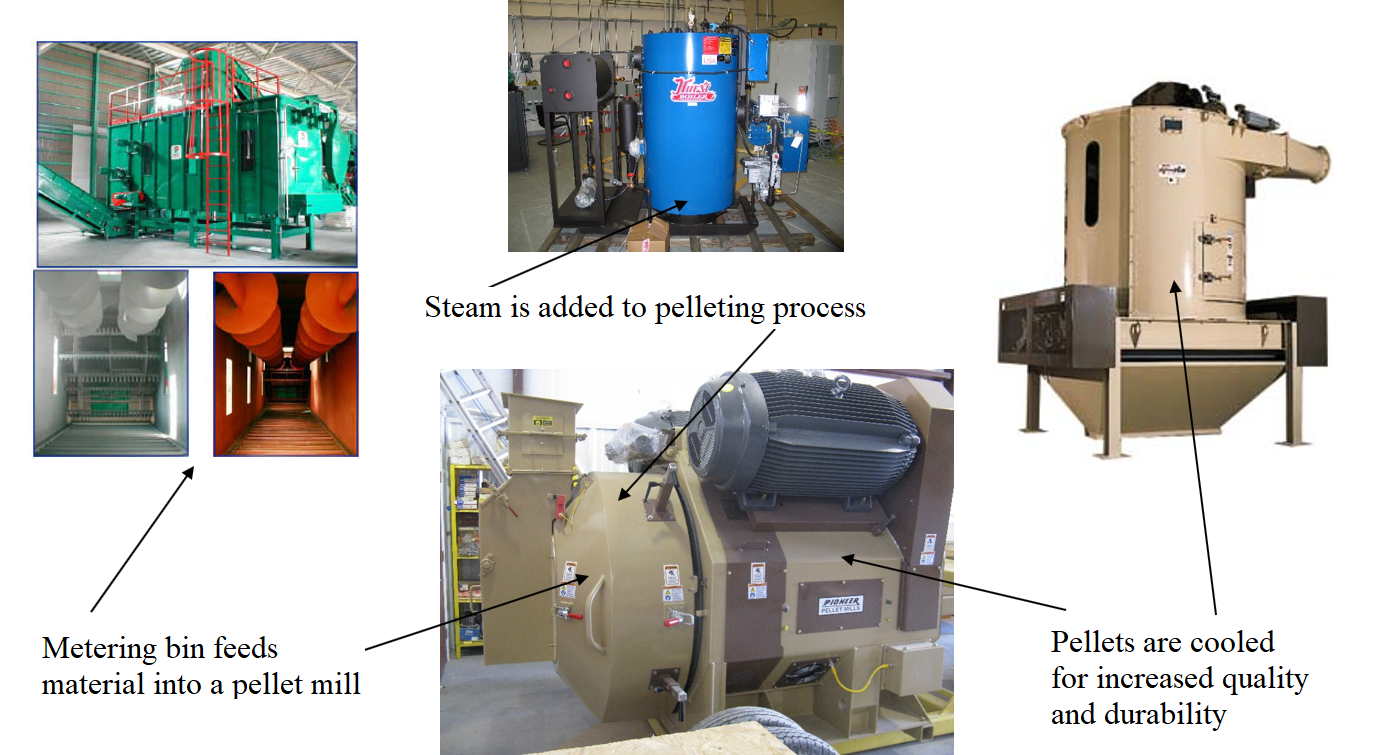}
	\caption{Pelleting equipment}
	\label{fig:pelleting}
\end{figure}

\begin{itemize}
	\itemsep-0.5em
	\item This is a Bliss Pioneer Pellet Mill.
	\item Material Input: Conveyors transfer ground material (smaller than 6.35 mm) to a metering bin which allows for even feeding into a pellet mill.
	\item Densification: It converts ground material from a metering bin into a dense product; goal is a minimum of 0.641 dry Mg/cubic meter. Baseline system includes a pellet mill, but other processes may be substituted in the future (examples are cubes or briquettes).
	\item Metering bin: Ground material must be fed into the pellet mill in a controlled fashion. The metering bin stores a supply of feed material, so it can be evenly conveyed into the mill. Although planned to feed material into the pellet mill, it can be used with other modules as needed.
	\item Conditioning: Steam provides moisture and heat to aid the densification process; other additives may be needed for certain materials or subsequent applications.
	\item Cooler: Cooling of the densified product improves the quality and durability. A cyclone separator and blower are connected to collect any fines created from the densified product.
	\item Material output: Pellets (or other product forms) will be conveyed into super sacks, transport trucks, or other storage containers.
\end{itemize}

\subsection{Belt Conveyors}

\begin{figure}[H]
	\centering
	\includegraphics[width=0.9\textwidth]{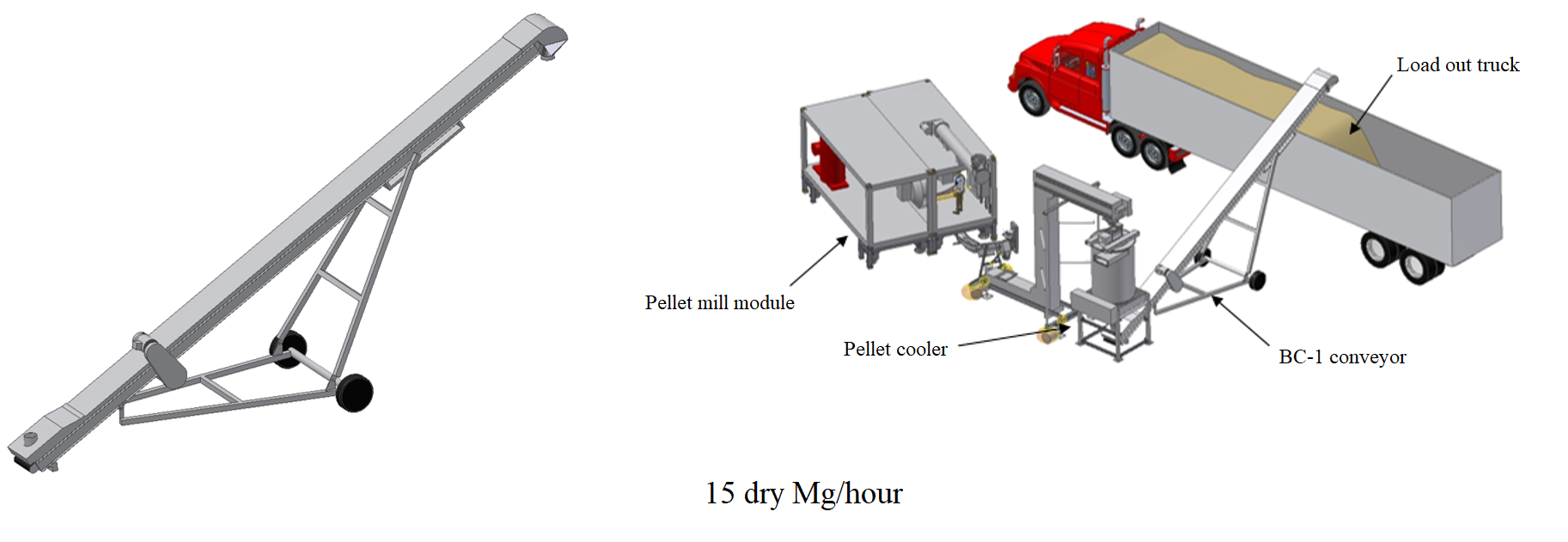}
	\caption{Belt conveyor}
	\label{fig:beltconveyor}
\end{figure}

\subsection{Drag Chain Conveyors}

\begin{figure}[H]
	\centering
	\includegraphics[width=0.9\linewidth]{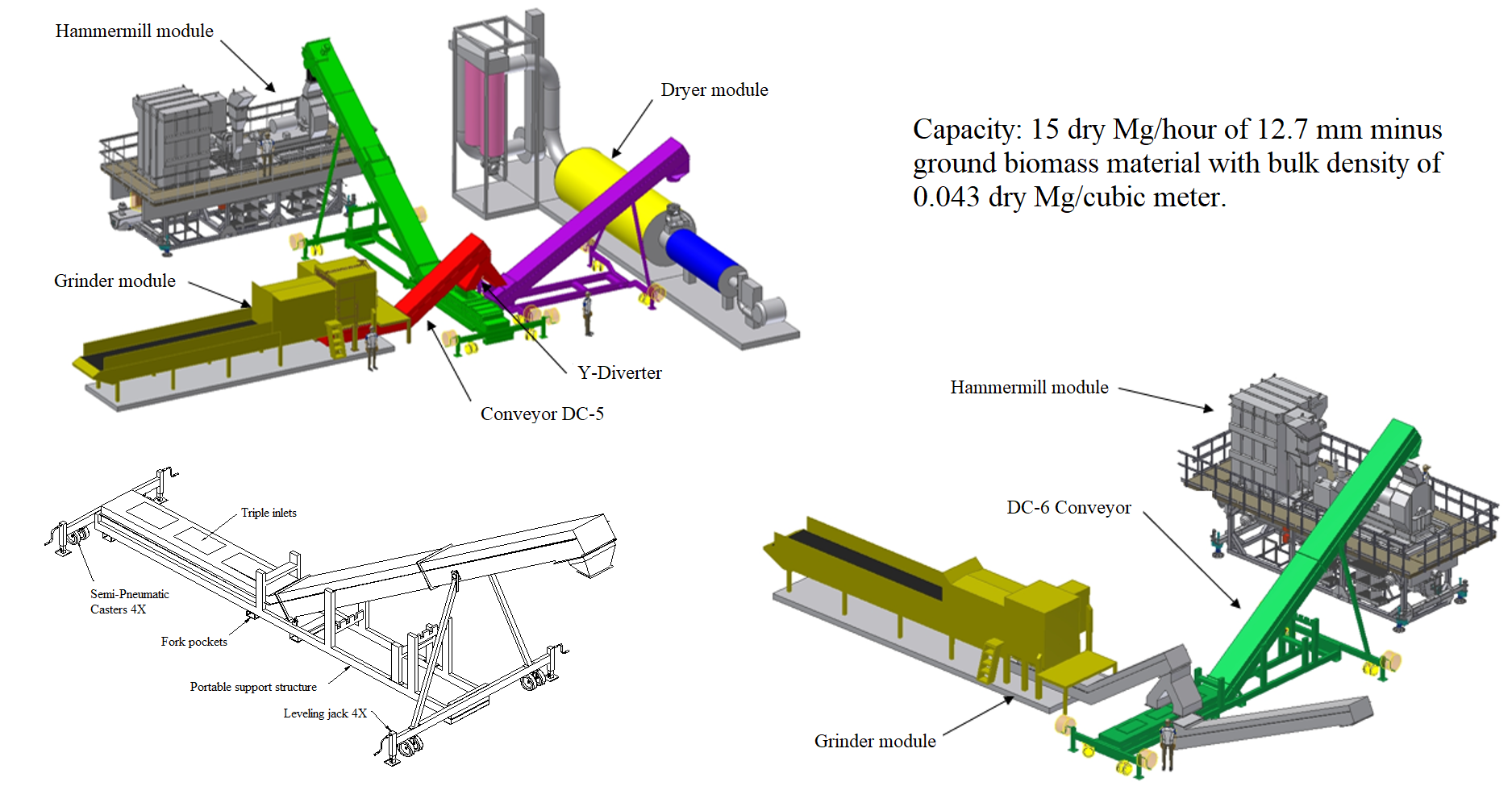}
	\caption{Drag chain conveyors}
	\label{fig:dragchain}
\end{figure}

\section{Input Data Utilized for Techno-economic Analysis}

\begin{table}[H]
	\footnotesize
	\centering
	\begin{threeparttable}
		\begin{tabular}{lccc}
			\toprule
			\multirow{2}{*}{\textbf{Equipment}} & \multirow{2}{*}{\textbf{Listed price}} & 
			\multirow{2}{*}{\textbf{kW}} &
			\textbf{Maximum operating} \\
			& & & \textbf{capacity (dry Mg/hour)} \\
			\midrule
			Bale conveyor & \$20,500 & 1.49 & 15 \\
			First-stage grinder & & 149.14	& 5 \\
			Drag chain conveyor 5 & \$49,733 & 2.24 & 15 \\
			Drag chain conveyor 6 &	\$65,949 & 3.73 & 15 \\
			Screw conveyor 1 & \$178,000 & 14.91 & 15 \\
			Second-stage grinder & & 111.86 & 5 \\ 
			Screw conveyor 2 & \$16,990 & 3.73 & 15 \\
			Screw conveyor 4 & \$58,824 & 3.73 & 15 \\
			Metering bin & \$186,264 & 2.98 & 24.07* \\
			Screw conveyor 5 & \$65,478 & 3.73 & 15 \\
			Pelleting equipment	& & 186.43 & 5 \\
			Drag chain conveyor 1 & \$31,146 & 1.49 & 15 \\
			Bucket elevator & \$50,000 & 1.49 & 15 \\	
			Counter flow & \$45,000 & 1.49 & 15 \\
			Belt conveyor & \$65,478 & 3.73 & 15 \\
			Storage bin	& \$186,264 & 2.98 & 5 \\
			Screw conveyor 6 & & & 15 \\	
			\bottomrule
		\end{tabular}
		\begin{tablenotes}
			\item * Operating capacity of metering bin is in cubic meter.
		\end{tablenotes}
	\end{threeparttable}
	\caption{Equipment specifications}
	\label{tab:eqspecs}
\end{table}

\begin{table}[H]
	\footnotesize
	\centering
	\begin{tabular}{lccclccc}
		\toprule
		\textbf{Unit cost of} & \multicolumn{3}{c}{\textbf{Moisture levels}} & \textbf{Unit cost of} & \multicolumn{3}{c}{\textbf{Moisture Levels}} \\
		\textbf{equipment (per hour)} & \textbf{High} & \textbf{Medium} & \textbf{Low} & \textbf{equipment (per hour)} & \textbf{High} & \textbf{Medium} & \textbf{Low} \\
		\midrule
		Bale conveyor & \$0.59 & \$0.59 & \$0.59 & First-stage grinder & \$33.04 & \$32.79 & \$31.79 \\
		Drag chain conveyor 5 & \$1.35 & \$1.35 & \$1.35 & Drag chain conveyor 6 & \$1.85 & \$1.85 & \$1.85 \\
		Screw conveyor 1 & \$1.38 & \$1.38 & \$1.38 & Second-stage grinder & \$17.18 & \$14.96 & \$14.83 \\
		Screw conveyor 2 & \$5.37 & \$5.37 & \$5.37 & Screw conveyor 4 & \$11.38 & \$11.38 & \$11.38 \\
		Metering bin & \$10.61 & \$10.61 & \$10.61 & Screw conveyor 5 & \$6.43 & \$6.43 & \$6.43 \\
		Pelleting equipment & \$17.65 & \$15.65 & \$14.96 & Drag chain conveyor 1 & \$1.72 & \$1.72 & \$1.72 \\
		Bucket elevator & \$2.70 & \$2.70 & \$2.70 & Counter flow cooler & \$0.78 & \$0.78 & \$0.54 \\
		Belt conveyor & \$2.51 & \$2.51 & \$2.51 & Storage bin & \$3.50 & \$3.50 & \$3.50 \\
		Screw conveyor 6 & \$11.38 & \$11.38 & \$11.38 & Separation unit & \$2.22 & \$2.22 & \$2.22 \\
		\bottomrule
	\end{tabular}
	\caption{Unit cost of operating equipment}
	\label{tab:costdata}
\end{table}

\section{Control Variables in HPC and BFFPC}

\begin{figure}[H]
	\centering
	\begin{subfigure}[b]{\linewidth}
		\centering
		\includegraphics[width=\textwidth]{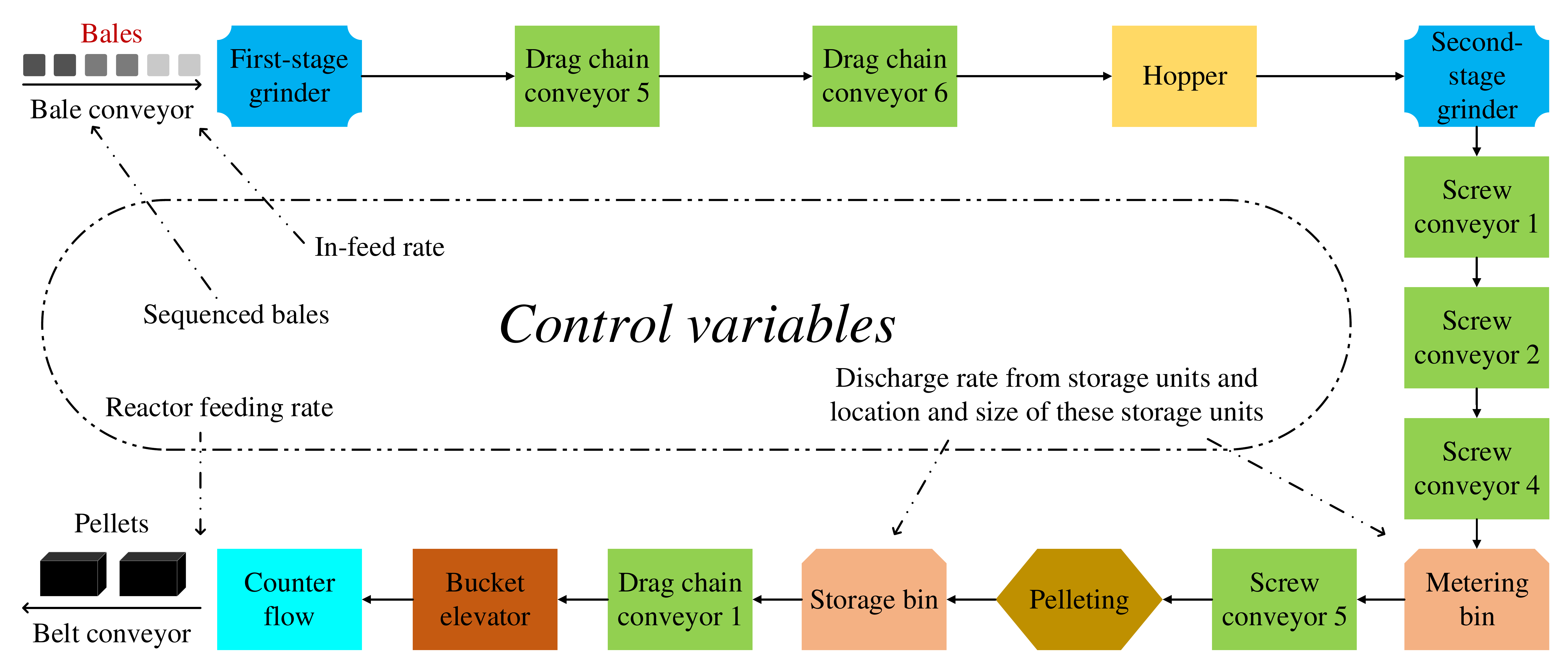}
		\caption{HPC}
		\label{fig:hpc_control}
		\vspace{2ex}
	\end{subfigure}
	\begin{subfigure}[b]{\linewidth}
		\centering
		\includegraphics[width=\textwidth]{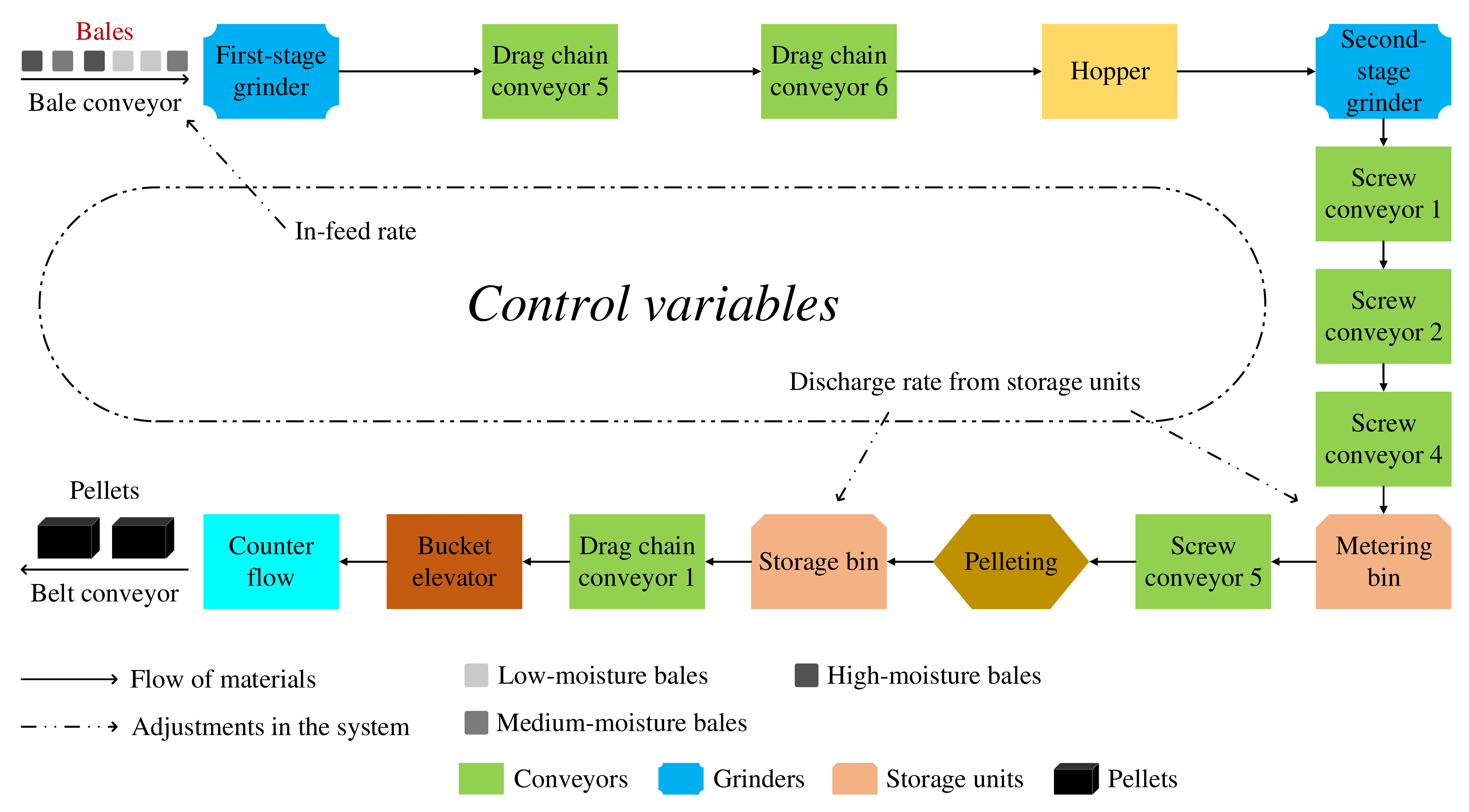}
		\caption{BFFPC}
		\label{fig:bffpc_control}
	\end{subfigure}
	\caption{Control variables at feedstock pre-processing system}
	\label{fig:control_var}
\end{figure}

\section{Commonalities and Differences of BFFPC and HPC}

\begin{table}[H]
	\footnotesize
	\centering
	\begin{tabularx}{\textwidth}{lXX}
		\toprule
		\textbf{Attribute} & \textbf{BFFPC} & \textbf{HPC} \\
		\midrule
		Control variables & Bale in-feed rate, discharge rates from storage units, discharge rate from pelleting & \textbf{Buffer size}, bale in-feed rate, discharge rates from storage units, discharge rate from pelleting \\
		\midrule
		State variables & \multicolumn{2}{l}{Amount of processed material in each location during each time step} \\
		\midrule
		Bale sequence & {Feeding pattern} is  \textbf{not} guided by {moisture level}. & {{{Feeding pattern}}} is guided by moisture level, maximization of throughput, and cost of operating system \\
		\midrule
		Feeding rate of reactor & Optimal control variables {\textbf{are not}} linked to minimization of feeding variability & Optimal control variables {\textbf{are linked}} to minimization of feeding variability \\
		\bottomrule
	\end{tabularx}
	\caption{Control approaches for feedstock pre-processing system and their attributes}
	\label{tab:control_summary}
\end{table}

\section{Identification of the Minimum Time to Process a Fixed Number of Bales}

\begin{figure}[H]
	\centering
	\includegraphics[width=\linewidth]{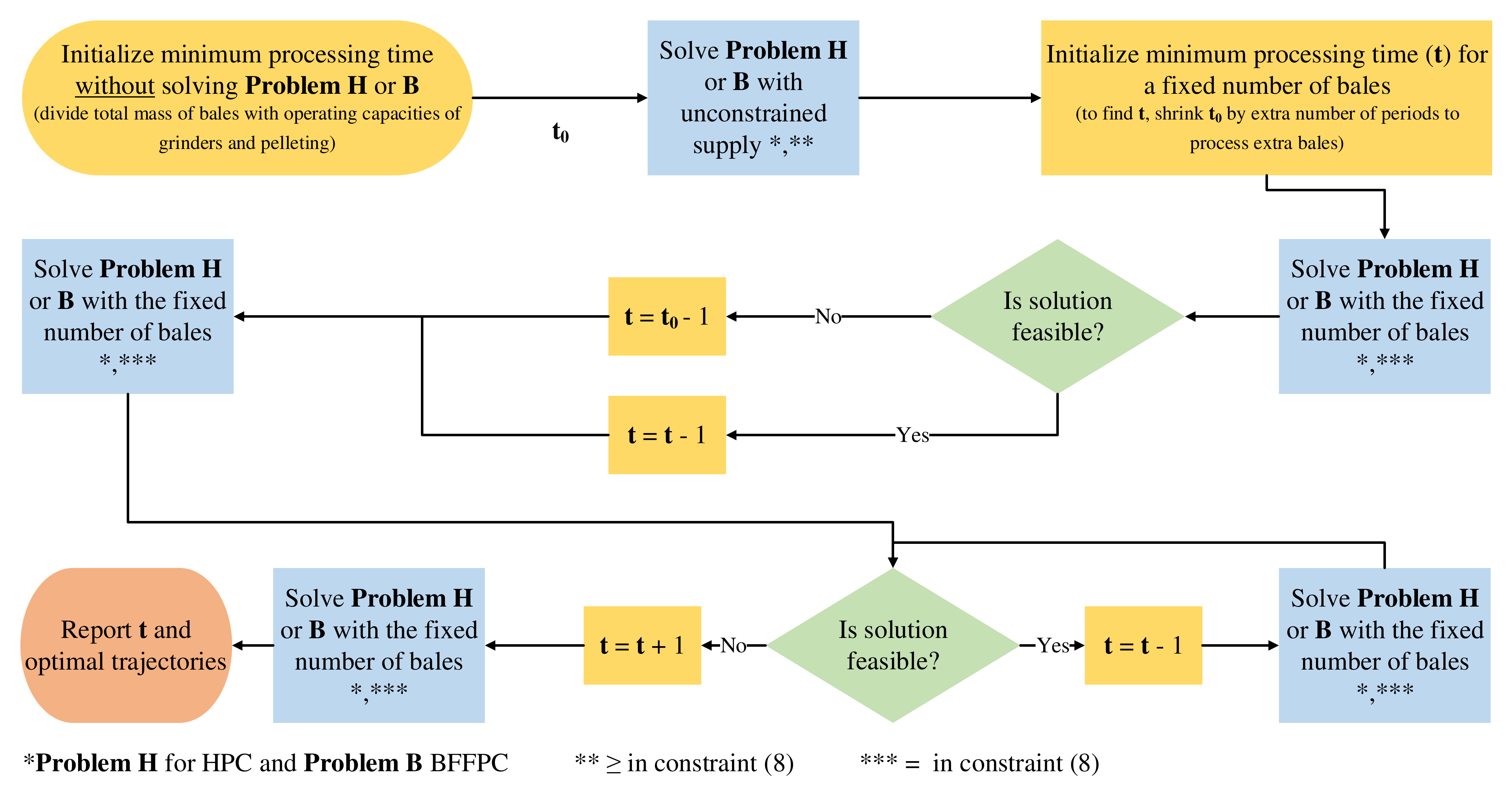}
	\caption{{An algorithm} to determine the minimum processing time and optimal trajectories to process a {fixed} number of bales}
	\label{fig:flowchart}
\end{figure}

The implementation of this algorithm is described below.

\subsection{Initialize minimum processing time} \label{sec:initmax}

In this step, we determine an initial modeling horizon ($t_{0}$) by dividing the total mass of bales (dry Mega gram (dry Mg)) with the system operating capacity (dry Mg/hour) for each moisture level. Then, we solve \textbf{Problem H} and \textbf{Problem B} to find the maximum number of bales that can be processed in HPC and BFFPC, respectively. These models are presented in Sections \ref{main-sec:hpcformulation} and \ref{main-sec:bffpcformulation} of the main text, respectively. In order to ensure infinite supply, we simply replace $=$ sign with $\geq$ sign in constraint \eqref{main-eq:allprocessed} of the main text. We note the solution of the optimization problem with unconstrained supply shown in Figure \ref{fig:flowchart} (blue box in the first row) is always feasible.

\subsection{Initialize minimum processing time for a fixed number of bales} \label{sec:initgiven}

This step initially finds the difference between the maximum and the required time to process a fixed number of bales. This difference is translated into extra time ($\delta$) assigned to processing bales. Next, it obtains an updated modeling horizon ($t$) by reducing $t_{0}$ by $\delta$. It lastly solves the optimization model with the updated modeling horizon. While solving the problem, we ensure that the feeding system processes exactly the same number of bales. That is, we solve the optimization problem with \eqref{main-eq:allprocessed} as an equality constraint.

Tasks explained in Section \ref{sec:initmax} and Section \ref{sec:initgiven} enable us to find a good starting modeling horizon.

\subsection{Improvement of initial minimum processing time for a fixed number of bales}

In the last step (shown as a loop in Figure \ref{fig:flowchart}), the {algorithm} iteratively reduces the modeling horizon by solving the optimization model. {The algorithm} terminates when the processing time identified is infeasible to process every bale.

Table \ref{tab:step1} presents a numerical example of the {algorithm} outlined in Section \ref{sec:initmax}. Assume that we solve the optimization model and obtain the following results: the maximum amount of biomass that can be processed in a fixed modeling horizon ($t_{0}$) is 25 dry Mg, 42 dry Mg, and 17 dry Mg for low-, medium-, and high-moisture materials. This table summarizes the calculations of the initial modeling horizon. Table \ref{tab:step2} summarizes the  calculations outlined in Section \ref{sec:initgiven}.

\begin{table}[H]
	\footnotesize
	\centering
	\begin{tabularx}{\textwidth}{lXXXX}
		\toprule
		\textbf{Moisture level} & \textbf{Number of bales to be processed} & \textbf{Total mass (dry Mg)} & \textbf{System capacity (dry Mg/hour)} & \textbf{Initial modeling horizon $t_{0}$ (hours)} \\
		\midrule
		Low & 60 & $60 \cdot 0.391 = 23.51$ & 4.76 & $23.51 / 4.76 = 4.94$ \\
		Medium & 100 & $100 \cdot 0.391 = 39.19$ & 2.80 & $39.19 / 2.80 = 13.98$ \\
		High & 40 & $40 \cdot 0.391 = 15.68$ & 1.59 & $15.68 / 1.59 = 9.87$ \\
		\bottomrule
	\end{tabularx}
	\caption{An example calculation for {the framework introduced} in Section \ref{sec:initmax}}
	\label{tab:step1}
\end{table}

\begin{table}[H]
	\footnotesize
	\centering
	\begin{tabularx}{\textwidth}{lXXXX}
		\toprule
		\textbf{Moisture level} & \textbf{Maximum amount of bales processed (dry Mg)} & \textbf{Extra amount of bales (dry Mg)} & \textbf{Extra time (hours)} & \textbf{Modeling horizon $t$ (hours)} \\
		\midrule
		Low & 25 & $25 - 23.51 = 1.49$ & $1.49 / 4.76 = 0.31$ & $4.94 - 0.31 = 4.62$ \\
		Medium & 42 & $42 - 39.19 = 2.81$ & $2.81 / 2.80 = 1.00$ & $13.98 - 1.00 = 12.98$ \\
		High & 17 & $17 - 15.68 = 1.32$ & $1.32 / 1.59 = 0.83$ & $9.87 - 0.83 = 9.04$ \\
		\bottomrule
	\end{tabularx}
	\caption{An example calculation for {the framework introduced} in Section \ref{sec:initgiven}}
	\label{tab:step2}
\end{table}